\documentclass[10.5pt]{article}
\usepackage{a4}

\usepackage{amsmath, amsfonts,amssymb,amscd,amsthm}
\usepackage{graphicx}

\usepackage{bm}
\usepackage{subfig}
\usepackage{graphicx}
\usepackage{pict2e}
\usepackage{color}
\usepackage[percent]{overpic}
\usepackage{verbatim}
\usepackage[normalem]{ulem}
\usepackage[colorlinks=true,linkcolor=black,urlcolor=black]{hyperref}
\usepackage[english]{babel}
\usepackage{algorithm, algorithmic}

\usepackage[title]{appendix}

\newtheorem{theorem}{Theorem}
\newtheorem{theorem*}{Theorem}

\newtheorem{remark}[theorem]{Remark}

\usepackage[utf8]{inputenc} 
\usepackage[T1]{fontenc}    
\usepackage{hyperref}       
\hypersetup{
    colorlinks,
    linkcolor={blue},
    citecolor={red},
    urlcolor={cyan}
}
\usepackage{url}            
\usepackage{booktabs}       
\usepackage{amsfonts}       
\usepackage{nicefrac}       
\usepackage{microtype}      
\usepackage{lipsum}
\usepackage{fancyhdr}       
\usepackage{graphicx}       
\usepackage{subfig}
\usepackage{wrapfig}
\graphicspath{{media/}}     

\newcommand{\MSC}[1]{%
\par\smallskip\noindent\textsc{Mathematics Subject Classification (MSC 2020):}~#1
}

\providecommand{\keywords}[1]{%
\par\smallskip\noindent\textsc{Keywords:}~#1
}

\usepackage{orcidlink}
\usepackage{tikz}
\usetikzlibrary{arrows.meta}
\usetikzlibrary{decorations.pathreplacing}

\newcommand*\xbar[1]{%
  \hbox{%
    \vbox{%
      \hrule height 0.5pt 
      \kern0.5ex
      \hbox{%
        \kern-0.1em
        \ensuremath{#1}%
        \kern-0.1em
    }%
    }%
  }%
} 
\definecolor{myblue}{rgb}{0.21, 0.34, 0.74}
\definecolor{myred}{rgb}{0.79, 0.0, 0.09}
\definecolor{mygreen}{rgb}{0, 0.32, 0}
\definecolor{mypurple}{rgb}{0.71, 0.36, 0.75}

\newcommand{\francesco}[1]{ {\colorbox{red}{\color{white} \textsf{FF}} \color{red}{#1}} }
\newcommand{\alessandra}[1]{ {\colorbox{blue}{\color{white} \textsf{AA}} \color{blue}{#1}} }

\begin{document}
\title{Finite elements for the space approximation\\ of a differential model for salts crystallization }

 \author{*${}^{\circ}$A. Aimi \orcidlink{0000-0002-4699-4261} , ${}^\dagger$${}^{\circ}$G. Bretti \orcidlink{0000-0001-5293-2115}, *${}^{\circ}$G. Di Credico \orcidlink{0000-0001-5339-0709}, $^\S$F. Freddi \orcidlink{0000-0003-0601-6022}, \\ *${}^{\circ}$C. Guardasoni \orcidlink{0000-0002-7054-8579} and ${}^{\ddagger\lozenge}$${}^{\circ}$M. Pezzella \orcidlink{0000-0002-1869-945X}\\ \\ \small *Dept.~of Mathematical Physical and Computer Sciences, University of Parma, Italy\\\small ${}^{\circ}$Members of the INdAM-GNCS Research Group, Italy\\
	\small$^\dagger\!\!$ National Research Council of Italy, Institute for Applied Mathematics ``M. Picone'', Rome branch, Italy\\
	\small$^\S$Department of Engineering and Architecture, University of Parma, Italy \\
	\small$^\ddagger\!\!$ National Research Council of Italy, Institute for Applied Mathematics ``M. Picone'', Naples branch, Italy\\
	\small$^\lozenge\!\!$ Dept. of Mathematics and Applications “R. Caccioppoli”, University of Naples Federico II,  Italy
}
 
\date{}

\maketitle 

\begin{abstract}
\noindent This article investigates a space-time differential model related to the degradation of stone artifacts caused by exposure to air and atmospheric agents, which specifically lead to the accumulation of salt crystals in the material. 
A numerical method based on finite-element space discretization and implicit-explicit time marching is proposed as an extension of the one-dimensional finite-difference framework introduced in \cite{BRACCIALE201721}. Within the same one-dimensional setting, a sensitivity analysis is performed, based on the techniques developed therein. They are also used as a comparison tool for the finite-element formulation, here introduced for more realistic simulations in higher space dimensions.
Considerations about stability will be provided, together with an experimental convergence analysis highlighting the performance of the proposed approach. Numerical results in two and three space dimensions, obtained by an efficient code implementation, will be presented and discussed.
\end{abstract}

\MSC{65M60 - 65M06 - 65Z05 - 76S05 - 93B35 }


\section{Introduction}
\label{intro}
Lapideous materials used in cultural heritage can be regarded as porous media, permeable to damaging agents transported by moisture which may lead to structural deterioration. When exposed to environmental factors, stones undergo weathering processes driven by water penetration, whether from meteoric precipitation or groundwater infiltration. Numerous studies have investigated the degradation of porous building materials in heritage structures (see, for instance, \cite{volume_mach, charola} and references therein). 

Degradation of monumental stones arises from a combination of physical-mechanical, chemical and biological processes. Physical-mechanical effects include salt crystallization \cite{doehne}, erosion \cite{hoke} and thermal expansion \cite{siegesmund}. Chemical reactions, such as carbonation \cite{ceseri, Freddi_2022}, sulfation \cite{sulfation, Freddi_2023}, oxidation \cite{reale} and hydration \cite{grondin, steiger}, alter the material composition, producing clay, crusts that swell, soluble salts and shrinkage. In addition, the growth of biological entities, including plants \cite{wang}, mosses, lichens and bacteria \cite{warscheid}, contributes to further deterioration. The interplay of these effects may induce microstructural changes in the material, sustaining an ongoing decay.

In this work, we focus on salt crystallization processes, which constitute one of the main causes of deterioration in porous materials. The phenomenon involves capillary penetration and internal reactions of salts dissolved in water flowing through the stone pores. Under certain conditions, the salts precipitate, forming crystals either on the surface (efflorescence) or within the pores (subflorescence). When crystallization occurs within the pores, it reduces the available pore volume, disrupts the liquid network and slows water transport, with consequent effects on the material deterioration.

In a context where preventive strategies play a critical role, several differential models have been proposed to characterize and quantify the degradation, alongside numerical methods specifically developed for its simulation. These models account for water transport in the porous medium, incorporating Darcy’s law, the changes in porosity due to crystal blooming and the effects of protective treatments aimed at mitigating the damage \cite{BRACCIALE201721, goidanich}. Existence and uniqueness results for the model in \cite{BRACCIALE201721}, both without and with boundary conditions, are provided in \cite{guarg, guarg2}. However, despite their modeling potential, the existing mathematical frameworks are based on a one-dimensional spatial representation. Although capable of describing vertical rise in porous materials, these models find limited applicability in realistic case studies. The first contribution of this work consists in the extension the model of \cite{BRACCIALE201721} to multidimensional domains, that provides a more detailed investigation of diffusion-driven phenomena. Furthermore, to address the challenges inherent in multidimensional simulations, a Finite Element Method (FEM) approach is employed. FEM provides a flexible and robust framework for the discretization of complex geometries and heterogeneous domains, ensuring accurate representation of the underlying physics. A comparison between FEM and traditional one-dimensional finite difference (FD) schemes, together with a discussion of their respective advantages and limitations, is reported in \cite{FDM-FEM}.

\subsection{General description of the model problem}
In this paper, we build upon
the mathematical model introduced in \cite{BRACCIALE201721}, describing the coupled processes of moisture transport, salt migration
and crystallization in porous materials exposed to saline water and atmospheric agents.
The formulation is motivated by standard laboratory imbibition--drying experiments performed on
stone- and brick-like specimens, where salt accumulation within the pore network leads to progressive
porosity reduction and material degradation.

The physical domain here considered 
is a bounded region $\Omega \subset \mathbb{R}^d$, with $d = 1,2,3$, representing
a prismatic porous specimen subjected to water absorption from the bottom surface and moisture exchange
with the surrounding air at the exposed boundaries.
The evolution of the system is studied over a finite time interval $[0,T]$, which may be decomposed
into an imbibition phase followed by a drying phase in order to reproduce the experimental protocol.

The model accounts for the interaction between a liquid phase carrying dissolved salt ions and a solid
crystalline phase precipitating within the pore space.
A key feature of the formulation is the explicit coupling between transport processes and microstructural
evolution induced by salt crystallization.

\subsection{State variables}

The unknowns of the problem are four space- and time-dependent fields defined for $(\textbf{x},t)\in\Omega \times [0,T]$:
\begin{itemize}
\item $\theta_l(\textbf{x},t)$, the volumetric fraction of liquid water, representing the portion of the pore
volume locally occupied by the liquid phase;

\item $c_i(\textbf{x},t)$, the concentration of dissolved salt ions transported by the liquid phase;

\item $c_s(\textbf{x},t)$, the concentration of crystallized salt deposited within the pore network;

\item $n(\textbf{x},t)$, the porosity of the material, defined as the fraction of the total volume available
for fluid transport.
\end{itemize}

The porosity $n$ is not treated as an independent field,
but is coupled to the crystallized salt
concentration through a constitutive relation accounting for pore clogging due to crystal growth.
As crystallization progresses, the available pore volume is reduced, which in turn affects the
transport properties of the medium.

\subsection{Governing mechanisms}

Moisture transport is described by a nonlinear diffusion process driven by capillary effects and
modulated by the evolving porosity.
The effective permeability of the medium depends on the saturation ratio $\theta_l/n$ through a
nonlinear function calibrated to reproduce experimentally observed absorption behavior in porous
materials.
This formulation allows the model to capture saturation-dependent mobility and sharp wetting fronts.

The transport of dissolved salt ions is governed by a convection--diffusion mechanism.
Advection is induced by the moisture flux, while molecular diffusion accounts for concentration
gradients within the liquid phase.
The ion balance is coupled to the crystallization kinetics, which acts as a sink term for the dissolved
salt and a source term for the solid phase.

Salt crystallization is modeled through a kinetic law depending on the local ion concentration and
the available pore space.
The precipitation of crystals reduces the porosity, which feeds back into the moisture and ion
transport equations, resulting in a strongly coupled nonlinear system.

\subsection{Experimental setting and scope}

Boundary conditions are chosen to reproduce standard imbibition--drying experiments.
During the imbibition phase, the bottom boundary of the specimen is assumed to be in contact with a
saline solution, prescribing saturated moisture and ion concentration.
At the exposed boundaries, moisture exchange with the external environment is modeled through
Robin-type conditions, while no-flux (homogeneous Neumann) conditions are imposed for the salt ions.
During the drying phase the boundary conditions are modified to simulate evaporation and the absence
of further salt supply.

The resulting model consists of a system of nonlinear, time-dependent partial differential equations
with strong coupling between transport, reaction, and microstructural evolution.
The formulation is posed in a general spatial dimension $d$, enabling one-dimensional simulations
for calibration and sensitivity analysis, as well as higher-dimensional computations for more
realistic geometries.

The aim of the model is to provide a predictive and numerically tractable framework for the analysis
of salt-induced degradation phenomena in porous materials, with particular emphasis on the
interplay between moisture dynamics, crystallization processes, and porosity variation.

The original contribution of the present work is therefore given by the following aspects: i)  a multi-dimensional extension of the 1D model in \cite{BRACCIALE201721} is provided; ii) a sensitivity analysis on model parameters is carried out in the 1D formulation;  iii) a stable and convergent FEM based solver is developed and implemented to simulate 2D and 3D geometries.\\ 

The paper is organized as follows: Section \ref{sec;the model problem} describes in detail the differential problem with its initial and boundary conditions, for both the imbibition and the drying phases; sensitivity analysis is developed in Section \ref{sec:Sensitivity}, while the description of the space-time discretization is written in Section \ref{sec: Time discretization and iterative weak problem}. Numerical results are presented and discussed in Section \ref{sec:numerical_results}, together with an experimental convergence analysis, while conclusions are collected in Section  \ref{sec:conclusions}. For the interested reader, the description of the calibration of the top-boundary exchange coefficient can be found in Appendix \ref{app:Kw}, while, to make the paper self-contained, the original 1D space-time finite differences scheme is briefly reported in Appendix \ref{app:FD}.


\section{The differential problem}
\label{sec;the model problem}


For the sake of simplicity,  the space domain of the model problem
will be represented by a $d$-dimensional prism $$\Omega=\left\lbrace \textbf{x}\in\mathbb{R}^d\:\vert\: (x_1,\ldots,x_{d-1})\in[-L/2,L/2]^{d-1},\: x_3\in[0,H] \right\rbrace,$$ being $L>0$ 
and $H>0$ its height. We adopt the convention that, if $d=1$, the physical domain reduces to the line interval $\Omega=[0,H]$ and the 1D space variable can be written equivalently as $x_3$ or $x$. The boundary $\partial \Omega=\Gamma$ is partitioned as the union of the sets $\Gamma_{\mathcal{T}}$, $\Gamma_{\mathcal{B}}$ and $\Gamma_{\mathcal{L}}$, where
$$\Gamma_{\mathcal{T}}=\left\lbrace \textbf{x}\in\Omega\:\vert\: x_3=H\right\rbrace,\quad\Gamma_{\mathcal{B}}=\left\lbrace \textbf{x}\in\Omega\:\vert\: x_3=0\right\rbrace,$$
are respectively the top base and the bottom base of the prism (reducing just to the two extremes of the interval $[0,H]$ for the 1D case), while the lateral edges of the prism are defined by $\Gamma_{\mathcal{L}}=\Gamma\setminus (\Gamma_{\mathcal{B}}\cup\Gamma_{\mathcal{T}})$, see Figure \ref{fig:3D experiment}.

\tikzset {_hbz17hkta/.code = {\pgfsetadditionalshadetransform{ \pgftransformshift{\pgfpoint{0 bp } { 0 bp }  }  \pgftransformrotate{-315 }  \pgftransformscale{2 }  }}}
\pgfdeclarehorizontalshading{_n354xv4zd}{150bp}{rgb(0bp)=(1,1,1);
rgb(37.5bp)=(1,1,1);
rgb(62.5bp)=(0.82,0.92,0.98);
rgb(100bp)=(0.82,0.92,0.98)}

\tikzset {_ab1ay2uaq/.code = {\pgfsetadditionalshadetransform{ \pgftransformshift{\pgfpoint{10 bp } { 10 bp }  }  \pgftransformrotate{-315 }  \pgftransformscale{2 }  }}}
\pgfdeclarehorizontalshading{_3fzu5d5ir}{150bp}{rgb(0bp)=(1,1,1);
rgb(37.5bp)=(1,1,1);
rgb(62.5bp)=(0.82,0.92,0.98);
rgb(100bp)=(0.82,0.92,0.98)}
\tikzset{every picture/.style={line width=0.75pt}} 
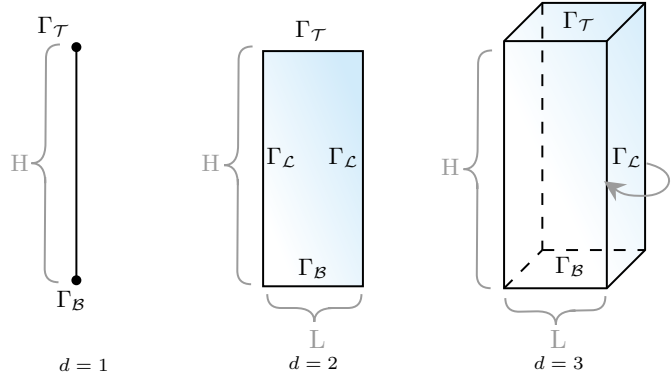
\begin{figure}
\centering
\begin{tikzpicture}[x=0.75pt,y=0.75pt,yscale=-1,xscale=1]
\draw [color={rgb, 255:red, 0; green, 0; blue, 0 }  ,draw opacity=1 ][line width=0.75]    (186.91,62) -- (186.91,179.67) ;
\draw [shift={(186.91,179.67)}, rotate = 90] [color={rgb, 255:red, 0; green, 0; blue, 0 }  ,draw opacity=1 ][fill={rgb, 255:red, 0; green, 0; blue, 0 }  ,fill opacity=1 ][line width=0.75]      (0, 0) circle [x radius= 2.01, y radius= 2.01]   ;
\draw [shift={(186.91,62)}, rotate = 90] [color={rgb, 255:red, 0; green, 0; blue, 0 }  ,draw opacity=1 ][fill={rgb, 255:red, 0; green, 0; blue, 0 }  ,fill opacity=1 ][line width=0.75]      (0, 0) circle [x radius= 2.01, y radius= 2.01]   ;
\path  [shading=_n354xv4zd,_hbz17hkta] (331.42,64) -- (331.42,182.71) -- (281.08,182.71) -- (281.08,64) -- cycle ; 
 \draw  [color={rgb, 255:red, 0; green, 0; blue, 0 }  ,draw opacity=1 ][line width=0.75]  (331.42,64) -- (331.42,182.71) -- (281.08,182.71) -- (281.08,64) -- cycle ; 
\path  [shading=_3fzu5d5ir,_ab1ay2uaq] (402.61,58.96) -- (421.91,39.67) -- (473.73,39.67) -- (473.73,164.42) -- (454.43,183.71) -- (402.61,183.71) -- cycle ; 
 \draw  [color={rgb, 255:red, 0; green, 0; blue, 0 }  ,draw opacity=1 ][line width=0.75]  (402.61,58.96) -- (421.91,39.67) -- (473.73,39.67) -- (473.73,164.42) -- (454.43,183.71) -- (402.61,183.71) -- cycle ; 
 \draw  [color={rgb, 255:red, 0; green, 0; blue, 0 }  ,draw opacity=1 ][line width=0.75]  (473.73,39.67) -- (454.43,58.96) -- (402.61,58.96) ; \draw  [color={rgb, 255:red, 0; green, 0; blue, 0 }  ,draw opacity=1 ][line width=0.75]  (454.43,58.96) -- (454.43,183.71) ;
\draw  [color={rgb, 255:red, 155; green, 155; blue, 155 }  ,draw opacity=1 ] (178.62,60.67) .. controls (173.95,60.64) and (171.61,62.95) .. (171.58,67.62) -- (171.32,110.62) .. controls (171.28,117.29) and (168.93,120.61) .. (164.26,120.58) .. controls (168.93,120.61) and (171.24,123.95) .. (171.19,130.62)(171.21,127.62) -- (170.93,173.62) .. controls (170.9,178.29) and (173.22,180.64) .. (177.89,180.67) ;
\draw  [color={rgb, 255:red, 155; green, 155; blue, 155 }  ,draw opacity=1 ] (275.76,61.67) .. controls (271.09,61.64) and (268.75,63.95) .. (268.72,68.62) -- (268.46,111.62) .. controls (268.42,118.29) and (266.07,121.61) .. (261.4,121.58) .. controls (266.07,121.61) and (268.38,124.95) .. (268.34,131.62)(268.35,128.62) -- (268.07,174.62) .. controls (268.04,179.29) and (270.36,181.64) .. (275.03,181.67) ;
\draw  [color={rgb, 255:red, 155; green, 155; blue, 155 }  ,draw opacity=1 ] (396.55,63.67) .. controls (391.88,63.64) and (389.54,65.95) .. (389.51,70.62) -- (389.25,113.62) .. controls (389.21,120.29) and (386.86,123.61) .. (382.19,123.58) .. controls (386.86,123.61) and (389.17,126.95) .. (389.13,133.62)(389.15,130.62) -- (388.86,176.62) .. controls (388.83,181.29) and (391.15,183.64) .. (395.82,183.67) ;
\draw  [color={rgb, 255:red, 155; green, 155; blue, 155 }  ,draw opacity=1 ] (283.18,186.67) .. controls (283.18,191.34) and (285.51,193.67) .. (290.18,193.67) -- (296.93,193.67) .. controls (303.6,193.67) and (306.93,196) .. (306.93,200.67) .. controls (306.93,196) and (310.26,193.67) .. (316.93,193.67)(313.93,193.67) -- (323.67,193.67) .. controls (328.34,193.67) and (330.67,191.34) .. (330.67,186.67) ;
\draw  [color={rgb, 255:red, 155; green, 155; blue, 155 }  ,draw opacity=1 ] (403.23,188.67) .. controls (403.32,193.34) and (405.7,195.62) .. (410.37,195.53) -- (418.6,195.36) .. controls (425.27,195.23) and (428.65,197.49) .. (428.74,202.16) .. controls (428.65,197.49) and (431.93,195.09) .. (438.6,194.96)(435.6,195.02) -- (446.83,194.8) .. controls (451.5,194.71) and (453.78,192.33) .. (453.69,187.66) ;
\draw [color={rgb, 255:red, 155; green, 155; blue, 155 }  ,draw opacity=1 ]   (474.47,120.67) .. controls (492.41,124.93) and (486.75,136.8) .. (473.72,137.04) .. controls (462.12,137.25) and (461.34,135.59) .. (455.41,131.36) ;
\draw [shift={(452.95,129.67)}, rotate = 33.29] [fill={rgb, 255:red, 155; green, 155; blue, 155 }  ,fill opacity=1 ][line width=0.08]  [draw opacity=0] (10.72,-5.15) -- (0,0) -- (10.72,5.15) -- (7.12,0) -- cycle    ;
\draw  [dash pattern={on 4.5pt off 4.5pt}]  (421.91,164.42) -- (473.73,164.42) ;
\draw  [dash pattern={on 4.5pt off 4.5pt}]  (402.61,183.71) -- (421.91,164.42) ;
\draw  [dash pattern={on 4.5pt off 4.5pt}]  (421.91,39.67) -- (421.91,164.42) ;
\draw (152.19,115) node [anchor=north west][inner sep=0.75pt]  [font=\small,color={rgb, 255:red, 155; green, 155; blue, 155 }  ,opacity=1 ] [align=left] {H};
\draw (248.59,115) node [anchor=north west][inner sep=0.75pt]  [font=\small,color={rgb, 255:red, 155; green, 155; blue, 155 }  ,opacity=1 ] [align=left] {H};
\draw (369.38,117) node [anchor=north west][inner sep=0.75pt]  [font=\small,color={rgb, 255:red, 155; green, 155; blue, 155 }  ,opacity=1 ] [align=left] {H};
\draw (302.79,203) node [anchor=north west][inner sep=0.75pt]  [color={rgb, 255:red, 155; green, 155; blue, 155 }  ,opacity=1 ] [align=left] {L};
\draw (424.33,203) node [anchor=north west][inner sep=0.75pt]  [color={rgb, 255:red, 155; green, 155; blue, 155 }  ,opacity=1 ] [align=left] {L};
\draw (165.78,45.4) node [anchor=north west][inner sep=0.75pt]  [font=\small]  {$\Gamma _{\mathcal{T}}$};
\draw (295.83,47.4) node [anchor=north west][inner sep=0.75pt]  [font=\small]  {$\Gamma _{\mathcal{T}}$};
\draw (430.78,42.4) node [anchor=north west][inner sep=0.75pt]  [font=\small]  {$\Gamma _{\mathcal{T}}$};
\draw (175.3,183.4) node [anchor=north west][inner sep=0.75pt]  [font=\small]  {$\Gamma _{\mathcal{B}}$};
\draw (297.34,167.4) node [anchor=north west][inner sep=0.75pt]  [font=\small]  {$\Gamma _{\mathcal{B}}$};
\draw (426.87,166.4) node [anchor=north west][inner sep=0.75pt]  [font=\small]  {$\Gamma _{\mathcal{B}}$};
\draw (312.51,111.4) node [anchor=north west][inner sep=0.75pt]  [font=\small]  {$\Gamma _{\mathcal{L}}$};
\draw (281.34,111.4) node [anchor=north west][inner sep=0.75pt]  [font=\small]  {$\Gamma _{\mathcal{L}}$};
\draw (455.56,111.4) node [anchor=north west][inner sep=0.75pt]  [font=\small]  {$\Gamma _{\mathcal{L}}$};
\draw (176.62,217.4) node [anchor=north west][inner sep=0.75pt]  [font=\scriptsize]  {$d=1$};
\draw (292.31,216.4) node [anchor=north west][inner sep=0.75pt]  [font=\scriptsize]  {$d=2$};
\draw (416.07,216.4) node [anchor=north west][inner sep=0.75pt]  [font=\scriptsize]  {$d=3$};
\end{tikzpicture}
\caption{Geometrical representation of the domain $\Omega$ for $d=1,2,3$, together with a schematic representation of the partition set on the boundary $\Gamma$.}
\label{fig:3D experiment}
\end{figure}



\subsection{Water absorption phase}
In this phase of the experiment, 
considered in the time domain $[0,T]$, the lower base $\Gamma_{\mathcal{B}}$ is supposed to be in contact with a salty water solution, seeping vertically along the height of the prism and depositing salt crystals in its volume. At the top base $\Gamma_{\mathcal{T}}$, interactions with the external environment are allowed, especially the drying action of the air. Imbibition from the bottom base induces the following sets of initial conditions, given for all dimensions $d=1,2,3$:
\begin{equation}\label{initial conditions}
    \begin{array}{cc}
         \theta_l(\textbf{x},0)=\begin{cases}
                     n_0\quad\textrm{for}\quad\textbf{x}\in\Gamma_{\mathcal{B}}\\[1pt]
                    \bar{\theta}_l\quad\textrm{for}\quad\textbf{x}\in\Omega\setminus\Gamma_{\mathcal{B}}
                                \end{cases}, &
        c_i(\textbf{x},0)=\begin{cases}
                    \bar{c_i}\quad\textrm{for}\quad\textbf{x}\in\Gamma_{\mathcal{B}}\\[1pt]
                    0\quad\textrm{for}\quad\textbf{x}\in\Omega\setminus\Gamma_{\mathcal{B}}
                                \end{cases},\\[12pt]
         c_s(\textbf{x},0)=0\quad\textrm{for}\quad\textbf{x}\in\Omega ,&
         
         n(\textbf{x},0)=n_0\quad\textrm{for}\quad\textbf{x}\in\Omega,
    \end{array}
\end{equation}
and the following Dirichlet constraints for $\theta_l$ and $c_i$ on $\Gamma_{\mathcal{B}}\times[0,T]:$
 \begin{equation}\label{dirichlet bc}
    \left\lbrace\begin{array}{cc}
    \theta_l(\textbf{x},t)=n_0&\textrm{for}\:(\textbf{x},t)\in\Gamma_{\mathcal{B}}\times[0,T]\\[2pt]
    c_i(\textbf{x},t)=\bar{c_i}&\textrm{for}\:(\textbf{x},t)\in\Gamma_{\mathcal{B}}\times[0,T]
    \end{array}\right..
\end{equation}
The positive parameters $n_0$ and $\bar{\theta}_l$ in \eqref{initial conditions} indicate, respectively, the porosity of the material and the moisture content in ambient air, while $\bar{c_i}$ is the sodium sulfate concentration in water.

Interaction with the air at the top of the brick is included by imposing Robin and Neumann conditions, respectively, for $\theta_l$ and  $c_i$:
\begin{equation}\label{Neumann-Robin bc}
    \left\lbrace\begin{array}{lc}
    \nabla \theta_l(\textbf{x},t)\cdot \textbf{n}=K_w \left(\bar{\theta}_l-\theta_l(\textbf{x},t) \right)&\textrm{for}\:(\textbf{x},t)\in\Gamma_{\mathcal{T}}\times[0,T]\\[2pt]
    \nabla c_i(\textbf{x},t)\cdot\textbf{n}=0&\textrm{for}\:(\textbf{x},t)\in \Gamma_{\mathcal{T}}\times[0,T]\\[2pt]
    \end{array}\right.,
\end{equation}
being $\textbf{n}:=\textbf{n}(\textbf{x})\in\mathbb{R}^d$ the normal outward pointing vector at the point $\textbf{x}\in\Gamma$. The loss of humidity at the top of the brick is addressed by the Robin constraint imposed on $\theta_l$ over $\Gamma_{\mathcal{T}}\times[0,T]$, with the positive constant $K_w$ controlling the intensity of the exchange with the external environment. This boundary treatment differs from that adopted in \cite{BRACCIALE201721}, where moisture exchange is described through an implicit nonlinear relation. Although physically detailed, that formulation is computationally demanding and does not naturally extend to multidimensional settings. Conversely, the choice in \eqref{Neumann-Robin bc} ensures well-posedness and straightforward applicability in higher dimensions. Consistency with the original model is retained by calibrating $K_w$ to reproduce, in the one-dimensional case, the same global dynamical behavior observed in \cite{BRACCIALE201721} for the two-phase imbibition–drying experiment. Details are provided in Appendix~\ref{app:Kw}.

Exclusively for $d>1$, the presence of the lateral sides $\Gamma_{\mathcal{L}}$ compels the definition of null flux constraints for quantities $\theta_l$ and $c_i$:
\begin{equation}\label{null lateral fluxes}
    \left\lbrace\begin{array}{l}
    \nabla \theta_l(\textbf{x},t)\cdot \textbf{n}=0\quad\textrm{for}\quad(\textbf{x},t)\in\Gamma_{\mathcal{L}}\times[0,T]\\[2pt]
    \nabla c_i(\textbf{x},t)\cdot\textbf{n}=0\quad\textrm{for}\quad(\textbf{x},t)\in \Gamma_{\mathcal{L}}\times[0,T]\\[2pt]
    \end{array}\right.,
\end{equation}
simulating the isolation of this portion of the boundary from the external environment, a condition obtainable in laboratory covering the porous prism by a layer of tape. Figure \ref{fig:Both_Phases} (left) provides a schematic three-dimensional representation of the imbibition phase.

\subsection{Drying phase}
The final quantities $\theta_l^T(\textbf{x}):=\theta_l(\textbf{x},T)$, $c_i^T(\textbf{x}):=c_i(\textbf{x},T)$, $c_s^T(\textbf{x}):=c_s(\textbf{x},T)$ and $n^T(\textbf{x}):=n(\textbf{x},T)$  obtained by the previous phase are used as initial conditions for the subsequent drying phase, whose time domain will be denoted by $[0,\widetilde{T}]$.  In fact, to be significant from the experimental point of view, this requires to investigate the dynamics starting from the spatial domain $\Omega$ saturated by humidity:

\begin{equation}\label{initial conditions drying}
\theta_l(\textbf{x},0)=\theta_l^T(\textbf{x}),\quad c_i(\textbf{x},0)=c_i^T(\textbf{x}),\quad c_s(\textbf{x},0)=c_s^T(\textbf{x}),\quad n(\textbf{x},0)=n^T(\textbf{x})\quad\textrm{for}\quad\textbf{x}\in\Omega.
\end{equation}

The simulation of the drying test requires a modification of the boundary conditions. In fact, we prescribe null flux for $c_i$ on the entire boundary, namely
\begin{equation}\label{bc_ci_drying}
    \nabla c_i(\textbf{x},t)\cdot\textbf{n}=0\quad\textrm{for}\quad(\textbf{x},t)\in\Gamma\times[0,\widetilde{T}].
    \end{equation}
For $\theta_l$ at the boundary, we prescribe 
\begin{equation}\label{conditions on theta test 1}
\left\lbrace\begin{array}{l}
    \theta_l(\textbf{x},t)=0\quad\textrm{for}\quad(\textbf{x},t)\in\Gamma_{\mathcal{B}}\cup\Gamma_\mathcal{T}\times[0,\widetilde{T}]\\[2pt]
    \nabla  \theta_l(\textbf{x},t)\cdot\textbf{n}=0\quad\textrm{for}\quad(\textbf{x},t)\in\Gamma_{\mathcal{L}}\times[0,\widetilde{T}]
    \end{array}\right..
\end{equation}

The Dirichlet constraints in \eqref{conditions on theta test 1} for $\theta$ at the top and bottom of the brick simulate the absence of humidity on the surface required to dry the brick.
We refer to Figure \ref{fig:Both_Phases} (right) for a schematic three-dimensional representation of the drying phase.

\tikzset {_8wrp6swyi/.code = {\pgfsetadditionalshadetransform{ \pgftransformshift{\pgfpoint{10 bp } { 10 bp }  }  \pgftransformrotate{-315 }  \pgftransformscale{2 }  }}}
\pgfdeclarehorizontalshading{_c0kx2v5ng}{150bp}{rgb(0bp)=(1,1,1);
rgb(37.5bp)=(1,1,1);
rgb(62.5bp)=(0.82,0.92,0.98);
rgb(100bp)=(0.82,0.92,0.98)}

\tikzset {_h9qa3oufd/.code = {\pgfsetadditionalshadetransform{ \pgftransformshift{\pgfpoint{10 bp } { 10 bp }  }  \pgftransformrotate{-315 }  \pgftransformscale{2 }  }}}
\pgfdeclarehorizontalshading{_7cmz4fy82}{150bp}{rgb(0bp)=(1,1,1);
rgb(37.5bp)=(1,1,1);
rgb(62.5bp)=(0.82,0.92,0.98);
rgb(100bp)=(0.82,0.92,0.98)}
\tikzset{every picture/.style={line width=0.75pt}} 
                
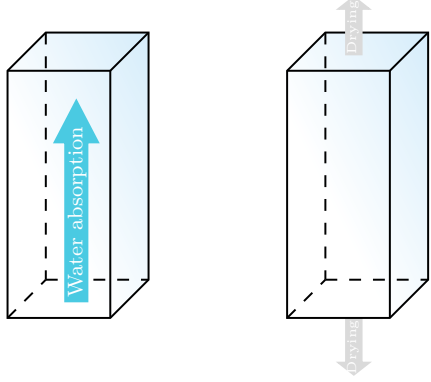
\begin{figure}
\centering
\begin{tikzpicture}[x=0.75pt,y=0.75pt,yscale=-1,xscale=1]
\draw  [color={rgb, 255:red, 218; green, 218; blue, 218 }  ,draw opacity=1 ][fill={rgb, 255:red, 218; green, 218; blue, 218 }  ,fill opacity=1 ] (355.35,188.33) -- (355.35,211.3) -- (351.57,211.3) -- (359.26,217.33) -- (366.94,211.3) -- (363.16,211.3) -- (363.16,188.33) -- cycle ;
\path  [shading=_c0kx2v5ng,_8wrp6swyi] (184.61,63.96) -- (203.91,44.67) -- (255.73,44.67) -- (255.73,169.42) -- (236.43,188.71) -- (184.61,188.71) -- cycle ; 
 \draw  [color={rgb, 255:red, 0; green, 0; blue, 0 }  ,draw opacity=1 ][line width=0.75]  (184.61,63.96) -- (203.91,44.67) -- (255.73,44.67) -- (255.73,169.42) -- (236.43,188.71) -- (184.61,188.71) -- cycle ; 
 \draw  [color={rgb, 255:red, 0; green, 0; blue, 0 }  ,draw opacity=1 ][line width=0.75]  (255.73,44.67) -- (236.43,63.96) -- (184.61,63.96) ; \draw  [color={rgb, 255:red, 0; green, 0; blue, 0 }  ,draw opacity=1 ][line width=0.75]  (236.43,63.96) -- (236.43,188.71) ;
\draw  [dash pattern={on 4.5pt off 4.5pt}]  (203.91,44.67) -- (203.91,169.42) ;
\draw  [dash pattern={on 4.5pt off 4.5pt}]  (184.61,188.71) -- (203.91,169.42) ;
\draw  [dash pattern={on 4.5pt off 4.5pt}]  (203.91,169.42) -- (255.73,169.42) ;
\draw  [color={rgb, 255:red, 74; green, 202; blue, 226 }  ,draw opacity=1 ][fill={rgb, 255:red, 74; green, 202; blue, 226 }  ,fill opacity=1 ] (213.71,180.33) -- (213.71,100.09) -- (208.32,100.09) -- (219.3,79) -- (230.28,100.09) -- (224.88,100.09) -- (224.88,180.33) -- cycle ;
\path  [shading=_7cmz4fy82,_h9qa3oufd] (325.61,63.96) -- (344.91,44.67) -- (396.73,44.67) -- (396.73,169.42) -- (377.43,188.71) -- (325.61,188.71) -- cycle ; 
 \draw  [color={rgb, 255:red, 0; green, 0; blue, 0 }  ,draw opacity=1 ][line width=0.75]  (325.61,63.96) -- (344.91,44.67) -- (396.73,44.67) -- (396.73,169.42) -- (377.43,188.71) -- (325.61,188.71) -- cycle ; 
 \draw  [color={rgb, 255:red, 0; green, 0; blue, 0 }  ,draw opacity=1 ][line width=0.75]  (396.73,44.67) -- (377.43,63.96) -- (325.61,63.96) ; \draw  [color={rgb, 255:red, 0; green, 0; blue, 0 }  ,draw opacity=1 ][line width=0.75]  (377.43,63.96) -- (377.43,188.71) ;
\draw  [dash pattern={on 4.5pt off 4.5pt}]  (344.91,169.42) -- (396.73,169.42) ;
\draw  [dash pattern={on 4.5pt off 4.5pt}]  (325.61,188.71) -- (344.91,169.42) ;
\draw  [dash pattern={on 4.5pt off 4.5pt}]  (344.91,44.67) -- (344.91,169.42) ;
\draw  [color={rgb, 255:red, 218; green, 218; blue, 218 }  ,draw opacity=1 ][fill={rgb, 255:red, 218; green, 218; blue, 218 }  ,fill opacity=1 ] (355.35,55.67) -- (355.35,32.7) -- (351.57,32.7) -- (359.26,26.67) -- (366.94,32.7) -- (363.16,32.7) -- (363.16,55.67) -- cycle ;
\draw (213.53,179.83) node [anchor=north west][inner sep=0.75pt]  [rotate=-270] [align=left] {\textcolor[rgb]{1,1,1}{{\footnotesize Water absorption}}};
\draw (354.35,57.00) node [anchor=north west][inner sep=0.75pt]  [font=\tiny,color={rgb, 255:red, 255; green, 255; blue, 255 }  ,opacity=1 ,rotate=-270] [align=left] {{\tiny Drying}};
\draw (354.35,218.3) node [anchor=north west][inner sep=0.75pt]  [font=\tiny,color={rgb, 255:red, 255; green, 255; blue, 255 }  ,opacity=1 ,rotate=-270] [align=left] {{\tiny Drying}};
\end{tikzpicture}
\caption{Schematic geometrical representation of the imbibition (left) and drying (right) phases.}
\label{fig:Both_Phases}
\end{figure}

\subsection{The governing equations}

The growth of salt particles in $\Omega$, both in water absorption and in drying dynamics, is defined by the following system of equations, taken from \cite{BRACCIALE201721}:

\begin{align}
        &\dfrac{\partial \theta_l}{\partial t}(\textbf{x},t)=\textrm{div}\!\left(\left(\dfrac{n(\textbf{x},t)}{n_0}\right)^2 \nabla B\!\left(\dfrac{\theta_l(\textbf{x},t)}{n(\textbf{x},t)}\right)\right),\label{cont theta}  \\[2pt]
        &\dfrac{\partial c_s}{\partial t}(\textbf{x},t)=K_s \, c_i(\textbf{x},t) \, (n(\textbf{x},t)-\theta_l(\textbf{x},t))^2 +\bar{K}\, (c_i(\textbf{x},t)-\bar{c})_+ \, \theta_l(\textbf{x},t)\label{cont cs}\\[2pt]
        &n(\textbf{x},t)=n_0-\gamma \, c_s(\textbf{x},t),  \hfill \quad (\textbf{x},t)\in \Omega\times[0,T],\label{cont n}\\[2pt]
        &\dfrac{\partial }{\partial t}\!\left(\theta_l(\textbf{x},t) \, c_i(\textbf{x},t)\right)=\textrm{div}\!\left(c_i(\textbf{x},t)\left(\dfrac{n(\textbf{x},t)}{n_0}\right)^2\nabla B\!\left(\dfrac{\theta_l(\textbf{x},t)}{n(\textbf{x},t)}\right)+D \, \theta_l(\textbf{x},t)\,  \nabla c_i(\textbf{x},t) \right)-\dfrac{\partial c_s}{\partial t}(\textbf{x},t),\label{cont ci}  \
\end{align}
where constants $K_s$ and $D$ are positive rates of crystallization and diffusivity, $\gamma$ indicates the crystal volume and $(\cdot)_+$ denotes the positive part of a function. The gradient operator $\nabla \left(\cdot\right)$ in \eqref{cont theta} and \eqref{cont ci} performs exclusively spatial derivatives and the function $B(s)$, defined as the following piece-wise third-degree polynomial function
\begin{equation}\label{eq:B_Forma}  
    B(s) =  
    \begin{cases}  
        0 & \text{if } s <a, \\  
        \dfrac{2}{3}c\left(\left(\dfrac{1-s}{1-a}\right)^2 \! (3a-1-2s) + (1-a)\right) \phantom{spazio} & \text{if } s \in [a,1],  \\
        \dfrac{2}{3}c(1-a) \phantom{spazio} & \text{if } s>1
    \end{cases}  
\end{equation}  
governs the fluid saturation in the brick. The parameters $a$ and $c$ in \eqref{eq:B_Forma} are positive constants related to the physical properties of the porous medium. The challenging identification of parameter sets of function $B(\cdot)$ that reproduces the capillary rise properties of specific porous materials is addressed in \cite{Data_Informed,Bretti2025}.\\
Let us finally remark that the equation \eqref{cont theta} describing the behavior of $\theta_l$ field is a nonlinear PDE.

In view of the approximation in space variable using FEM, let us rewrite \eqref{cont theta} and \eqref{cont ci} in the following compact notation 
\begin{align*}
&\dfrac{\partial \theta_l}{\partial t}=
       \textrm{div}\left[f(\theta_l,n)\nabla\theta_l\right]-\textrm{div}\left[F\left(\theta_l,n\right)\theta_l\right]\\     
       &\dfrac{\partial }{\partial t}\!\left(\theta_l\, c_i\right)=\textrm{div}\!\left[c_i\left( f(\theta_l,n)\nabla \theta_l - F(\theta_l,n) \theta_l \right)   \right]+D \,\textrm{div}\!\left[\theta_l\cdot\nabla c_i\right]-\dfrac{\partial c_s}{\partial t},
\end{align*}
being $f(\textbf{x},t)$ and $F(\textbf{x},t)$ the scalar and the vectorial operators
\begin{align*}
        &f\left(\theta_l(\textbf{x},t),n(\textbf{x},t)\right)=\dfrac{n(\textbf{x},t)}{n_0^2}\: B'\!\left(\dfrac{\theta_l(\textbf{x},t)}{n(\textbf{x},t)}\right),  \\[2pt]
        &F\left(\theta_l(\textbf{x},t),n(\textbf{x},t)\right)=\dfrac{1}{n_0^2}\: B'\!\left(\dfrac{\theta_l(\textbf{x},t)}{n(\textbf{x},t)}\right)\nabla n(\textbf{x},t).
\end{align*}
The exact expression of $B(s)$ derivative is 
\begin{equation}\label{eq:B_Forma_der}  
    B'(s) =  
    \begin{cases}  
        0 & \text{if } s <a\vee s>1, \\  
        \dfrac{4c}{(1-a)^2}(1-s)(s-a) \phantom{spazio} & \text{if } s \in [a,1].  
    \end{cases}  
\end{equation}  

We conclude the Section listing in Table \ref{tab:Old-Imbibition-Experiment} the values of the physical constants, together with their units of measurement, onto which the model problem \eqref{initial conditions}-\eqref{eq:B_Forma_der} depends.
\begin{table}[htbp]
 \caption{Parameters of the imbibition phenomenon simulation as detailed in \cite{Clarelli_Natalini_Nitsch}. The value of $K_w$ is determined as described in Appendix \ref{app:Kw}.
 }
  \centering
  \begin{tabular}{llll}
    \toprule
    & \textbf{Description}  & \textbf{Value}   & \textbf{Units} \\
    \midrule
    $n_0$     & Initial porosity of the material       & $2.8510 \cdot 10^{-1}$ & ---  \\
    $c$     & Physical property of the porous matrix       & $9.8073 \cdot 10^{-4}$ & cm$^2/$s  \\
    $a$     & Minimum saturation level for hydraulic continuity       & $2.1904 \cdot 10^{-1}$ & ---  \\
    $D$     & Na$_2$SO$_4$ diffusivity       & $1.2300 \cdot 10 ^{-5}$ & cm$^2/$s \\ 
    $\bar{\theta}_l$     & Moisture content of the ambient air       & $6.2540 \cdot 10^{-2}$ & g$/$cm$^3$ \\ 
    $\bar c_i$ &  Concentration in water of sodium sulphate & $9.9500 \cdot 10^{-2}$ & g$/$cm$^3$ \\ 
    $\gamma$ & Specific volume of crystals & $6.0000 \cdot 10^{-1}$ & cm$^3/$g \\
    $K_s$ & Crystallization rate coefficient & $4.1000 \cdot 10^{-5}$ & s$^{-1}$\\
    $K_w$ & Evaporation coefficient & $1.5\cdot 10^{-2}$ & cm$/$s \\
     $\bar{c}$ & Saturated concentration in water of sodium sulphate & $0.4399$ & g$/$cm$^3$  \\
      $\xbar K$ & Growth rate of hydrated crystals & $1.0000 \cdot 10^{-4}$ & s$^{-1}$ \\
    \bottomrule
  \end{tabular}
  \label{tab:Old-Imbibition-Experiment}
\end{table}

\section{Sensitivity analysis}\label{sec:Sensitivity}

In order to assess the robustness of the equations \eqref{cont theta}--\eqref{cont ci} with respect to the parameters that exert the strongest influence on the crystallization dynamics, we carried out a local sensitivity analysis following the approach of \cite{CdS}. Unlike classical one-at-a-time (OAT) investigations, which vary a single parameter while keeping others fixed, we explore local variations of multiple parameters simultaneously within prescribed bounds. In practice, this corresponds to probing a full neighborhood around the reference parameter set, rather than individual coordinate lines, thus capturing all relevant directional derivatives.

In this context, the 1D setting offers a computationally efficient framework. Specifically, we focus on the following three parameters 
\begin{itemize}
    \item $\gamma$, the specific volume of the crystals.  
    It determines the amount of pore volume locally occupied by the precipitated solid phase and directly influences the porosity reduction.

    \item $K_s$, the crystallization rate.  
    This coefficient governs the speed at which dissolved ions transform into solid crystals and therefore controls the temporal scale of salt accumulation.

    \item $K_{w}$, the exchange coefficient appearing in the Robin-type boundary condition \eqref{Neumann-Robin bc} at the top.  
\end{itemize}
The analysis evaluates the effect of perturbations of these parameters on the full evolution of the system, encompassing both the imbibition and drying phases. As sensitivity metrics, we consider the global quantities
\begin{equation}\label{eq:Sensitivity_Metrics}
    \mathcal{N}(\gamma,K_s,K_{w})= \frac{1}{H}\int_0^H n(x,T_{\mathrm{end}})\,dx,
    \qquad
    \mathcal{C}_s(\gamma,K_s,K_{w})= \frac{1}{H}\int_0^H c_s(x,T_{\mathrm{end}})\,dx,
\end{equation}
which represent, respectively, the average porosity and the average amount of crystallized salt at the final time $T_{\mathrm{end}} = T_{\mathrm{imb}} + T_{\mathrm{dry}} = 10\;\text{days} + 5\;\text{hours}.$ For both metrics, discretized by a fourth order Gregory quadrature rule (see, for instance, \cite{Gregory,Gregory2} and references therein), we evaluate the relative deviation with respect to the baseline parameters $(\gamma^0, K_s^0, K_{w}^0)$. The coefficients $\gamma^0=6\cdot 10^{-1} \ \textrm{cm}^3\textrm{/g}$ and $K_s^0=4.1\cdot 10^{-5} \ \textrm{s}$ coincide with the ones calibrated in \cite{BRACCIALE201721} (see Table \ref{tab:Old-Imbibition-Experiment}), whereas $K_{w}^0=K^*_w=1.5\cdot 10^{-2}$ is the optimal value identified with the analysis performed in Appendix~\ref{app:Kw}. 

\begin{figure}[htbp]
  \centering
    \includegraphics[width=0.9\linewidth]{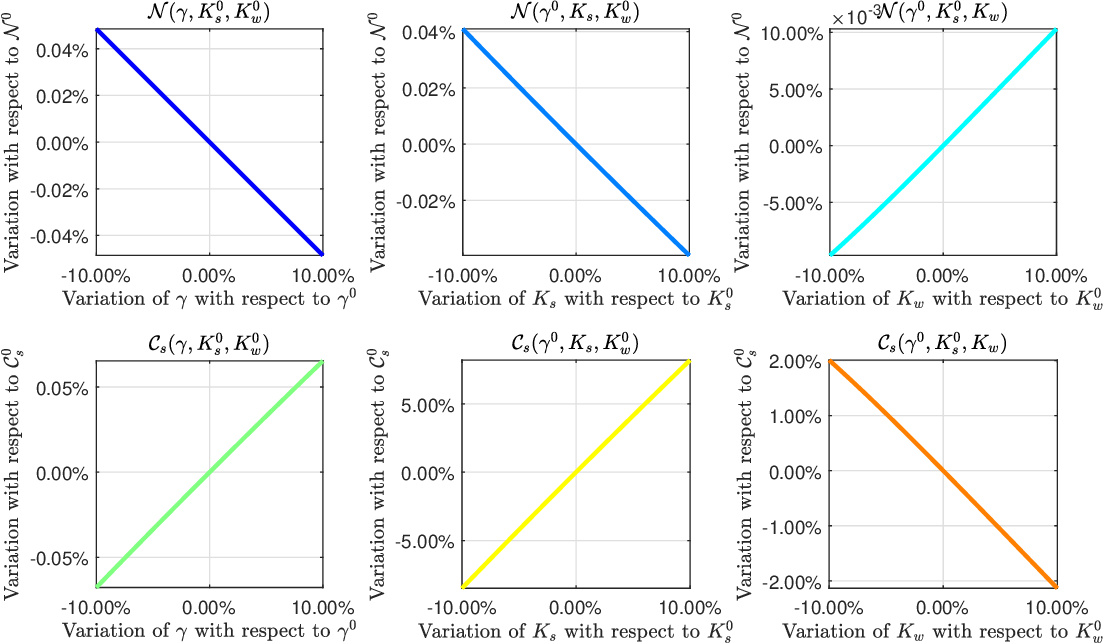} 
  \caption{Results of the one-factor-at-a-time (OAT) sensitivity analysis. 
  }
  \label{fig:OAT}
\end{figure}

To explore the parameter space, we account for increments and decrements of $1\%$ up to a maximum of $10\%$ for each parameter, which yields $21^3=9261$ parameter triplets. For each configuration, the complete two-phase simulation is performed up to $T_{\mathrm{end}}=882000\;\mathrm{s}$, using the finite difference scheme \eqref{eq:Num_method_1D} (see Appendix~\ref{app:FD}) with space and time steps $\Delta x=0.15$ cm and $\Delta t=0.25$ s. For each run, we compute the relative changes in $\mathcal{N}$ and $\mathcal{C}_s$ with respect to the reference values $\mathcal{N}^0=\mathcal{N}(\gamma^0,K_s^0,K_{w}^0)$ and $\mathcal{C}_s^0=\mathcal{C}_s^0(\gamma^0,K_s^0,K_{w}^0).$
\begin{figure}[htbp]
  \centering
    \includegraphics[width=0.9\linewidth]{media/N_3D_Sensitivity.png} \\[0.5cm]
    \includegraphics[width=0.9\linewidth]{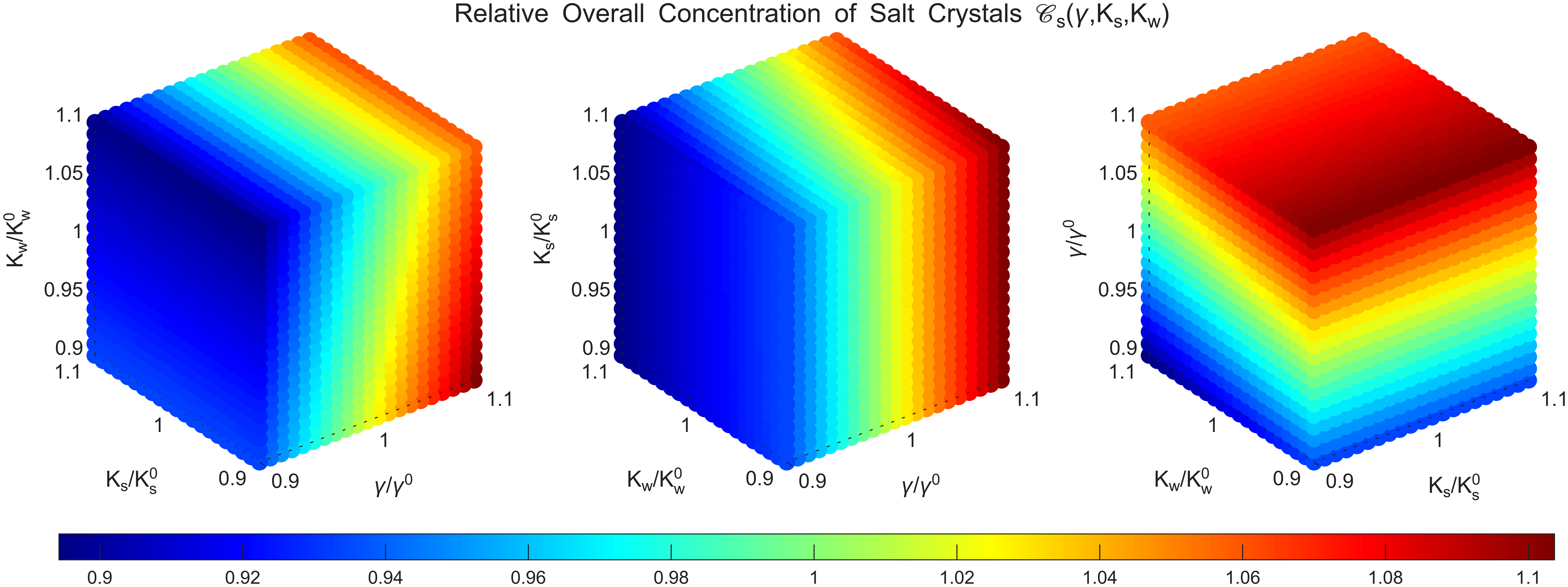}
  \caption{Results of the local sensitivity analysis. Relative changes, with respect to the baseline parameters, in the average porosity (top) and in the average crystallized salt (bottom) at the final time. Each plot illustrates the combined effect of perturbations in $\gamma$, $K_s$, and $K_w$. (Distinct colormaps and scales are used for the two panels).}
  \label{fig:3D_Sensitivity}
\end{figure}
Figure \ref{fig:3D_Sensitivity} reports the results of the local sensitivity analysis, performed as specified above. The model demonstrates strong robustness: even with parameter perturbations up to $\pm 10 \%$, the maximum deviation in $\mathcal{N}$ is $0.94 \%$, while for $\mathcal{C}_s$ it is $10.37 \%$. In both cases, the perturbations do not amplify in the system response. Furthermore, following the one-factor-at-a-time (OAT) approach \cite{OAT}, of which our methodology is an extension, we observe from Figure \ref{fig:OAT} that the water exchange coefficient $K_w$ has the largest influence on $\mathcal{N}$, whereas $\mathcal{C}_s$ is more sensitive to variations in the crystallization rate $K_s$. These results are in full agreement with the physical mechanisms governing the phenomenon.

\begin{remark}
   In line with the sensitivity analysis discussed above, the calibrated coefficients display stable behavior when applied within the full dynamics. This observation indicates that the parameter values identified in \cite{BRACCIALE201721}, despite being obtained in a different experimental and geometrical configuration, remain suitable for the multidimensional simulations conducted in this work.
\end{remark}

\section{Finite elements space discretization}\label{sec: Time discretization and iterative weak problem}
In this section we provide a hybrid elaboration of \eqref{cont theta}-\eqref{cont ci}. First of all, we introduce the time discretization: let $\Delta t$ be the positive time step size and consider the uniform mesh $t_k=k \Delta t,$ for $k=0,\ldots, N$, being $N$ such that the final time instant of analysis is given by $N \Delta t$. The four problem unknowns are then approximated by their collocation in time at the generic node $t_k$:  
\begin{equation}\label{temporal samples}
   \theta^{k}_l(\textbf{x})\approx \theta_l(\textbf{x}, t_k), \quad c_s^{k}(\textbf{x})\approx c_s(\textbf{x}, t_k), \quad  c_i^k(\textbf{x})\approx c_i(\textbf{x}, t_k), \quad n^{k}(\textbf{x})\approx n(\textbf{x}, t_k), \quad \ k=0,\ldots,N,
\end{equation}
allowing to approximate the involved time derivatives in terms of finite differences of the first order
\begin{equation}\label{discrete time derivatives}
\dfrac{\partial \theta_l}{\partial t}(\textbf{x},t_{k+1})\approx \frac{\theta^{k+1}_l(\textbf{x})-\theta^{k}_l(\textbf{x})}{\Delta t},\quad \dfrac{\partial c_s}{\partial t}(\textbf{x},t_{k+1})\approx \frac{c_s^{k+1}(\textbf{x})-c_s^{k}(\textbf{x})}{\Delta t},\quad  \dfrac{\partial c_i}{\partial t}(\textbf{x},t_{k+1})\approx\frac{c_i^{k+1}(\textbf{x})-c_i^{k}(\textbf{x})}{\Delta t}.
\end{equation}

It is understood that time approximations in \eqref{temporal samples} are subject to initial conditions
\eqref{initial conditions} for the imbibition phase and to 
\eqref{initial conditions drying} for the drying phase.\\

The space discretization starts writing the PDEs of the model problem, already discretized in time variable, in weak form. To do this, for the imbibition phase we consider test functions $v_{\theta_l},v_{c_i}\in H^1_{\Gamma_{\mathcal{B}}}(\Omega)=\left\lbrace v\in H^1(\Omega)\:\vert\: v_{\vert_{\Gamma_{\mathcal{B}}}}=0 \right\rbrace$ and we perform the projection of the PDEs \eqref{cont theta} and \eqref{cont ci} at each time step, obtaining for each $k=0,\ldots,N-1$ the following system of equations in the space unknowns $\theta_l^{k+1},c_i^{k+1},c_s^{k+1}$ and $n^{k+1}$, subjected to \eqref{dirichlet bc}-\eqref{Neumann-Robin bc}:

\begin{align}
        &\biggl\langle\frac{\theta^{k+1}_l-\theta^{k}_l}{\Delta t},v_{\theta_l}\biggr\rangle_{\Omega}=
        -\biggl\langle f(\theta^k_l,n^k)\nabla \theta^{k+1}_l,\nabla v_{\theta_l}\biggr\rangle_{\Omega}+\biggl\langle f(\theta^k_l,n^k) K_w \left(\bar{\theta}_l-\theta^{k+1}_l \right),v_{\theta_l} \biggr\rangle_{\Gamma_{\mathcal{T}}}\nonumber\\
        & \quad\quad\quad\quad\quad\quad\quad\quad\quad\quad\quad\quad\quad\quad\quad\quad\quad\quad-\biggl\langle\textrm{div}\left[F(\theta^k_l,n^k)\theta^{k+1}_l\right],v_{\theta_l}\biggr\rangle_{\Omega},\label{theta k} \\[4pt]
                &\frac{c_s^{k+1}-c_s^{k}}{\Delta t}=K_s \, c_i^{k} \, (n^{k}-\theta^{k}_l)^2+\bar{K}\, (c_i^k-\bar{c})_+ \, \theta^k_l\label{cs k},\\[4pt]
        &n^{k+1}=n_0-\gamma \, c_s^{k+1},\label{n k}\\[4pt]
   &\biggl\langle\frac{\theta^{k+1}_lc_i^{k+1}-\theta^{k}_lc_i^k}{\Delta t},v_{c_i}\biggr\rangle_{\Omega}=\biggl\langle {\rm div} \left\{c_i^{k+1}\left[ f(\theta^{k+1}_l,n^{k+1})\nabla \theta^{k+1}_l - \theta^{k+1}_l F(\theta^{k+1}_l,n^{k+1}) \right]\right\}, v_{c_i}\biggr\rangle_{\Omega}\nonumber\\[4pt]
      &\quad\quad\quad\quad\quad\quad\quad\quad\quad  -\biggl\langle D \theta^{k+1}_l\nabla c_i^{k+1} , \nabla v_{c_i} \biggr\rangle_{\Omega} - \biggl\langle \frac{c_s^{k+1}-c_s^{k}}{\Delta t},v_{c_i}\biggr\rangle_{\Omega}.\label{ci k}
\end{align}
The time iterations of the scheme starts from $\theta_l^0,c_i^0,c_s^0$ and $n^0$ given by the initial conditions \eqref{initial conditions}. 


For the drying phase, let us note that $v_{\theta_l}\in H^1_{\Gamma_{\mathcal{B}}\cup\Gamma_{\mathcal{{\cal T}}}}(\Omega)=\left\lbrace v\in H^1(\Omega)\:\vert\: v_{\vert_{\Gamma_{\mathcal{B}}\cup \Gamma_{\mathcal{{\cal T}}}}}=0 \right\rbrace$ and $v_{c_i}\in H^1(\Omega)$. Moreover, we remark that boundary conditions \eqref{conditions on theta test 1} induce a modification in the weak problem to be solved at each time step: in particular equation \eqref{theta k} simplifies in
\begin{align}
        &\biggl\langle\frac{\theta^{k+1}_l-\theta^{k}_l}{\Delta_t},v_{\theta_l}\biggr\rangle_{\Omega}=
        -\biggl\langle f(\theta^k_l,n^k)\nabla \theta^{k+1}_l,\nabla v_{\theta_l}\biggr\rangle_{\Omega}-\biggl\langle\textrm{div}\left[F(\theta^k_l,n^k)\theta^{k+1}_l\right],v_{\theta_l}\biggr\rangle_{\Omega}.\label{theta k drying test 1} 
\end{align}
In this phase, the time iterations start from $\theta_l^0,c_i^0,c_s^0$ and $n^0$ given by the initial conditions \eqref{initial conditions drying}.\\

In the above framework, the space discretization is performed introducing a standard conforming domain mesh, depending on the geometric parameter $h$, which is defined as the the maximum among triangles or tetrahedrons diameters, respectively, in 2D or 3D (and simplifies as the maximum length of the mesh elements in 1D). Finite dimensional spaces of piece-wise linear finite elements, related to the considered mesh, have been taken into account, to fully approximate the above equations. The hybrid first order FD in time and FEM in space has been set up within FEniCS computing platform \cite{FEniCS}.

\begin{remark}
Let us observe that the weak equation \eqref{theta k} presents an implicit-explicit discretization scheme in time for the $\theta_l$ unknown. Moreover, the fully discretized equation coming from \eqref{ci k} does not need the addition of a stabilizing term, as in the fourth equation of the finite differences scheme \eqref{eq:Num_method_1D} (see Appendix~\ref{app:FD}).  
\end{remark}

\section{Numerical results}
\label{sec:numerical_results}

\subsection{The 2D bar: experiment for long time imbibition phase}\label{The 2D bar: esperiment 1 for imbibition phase}
\begin{figure}[h!]
\centering
\includegraphics[scale=0.62]{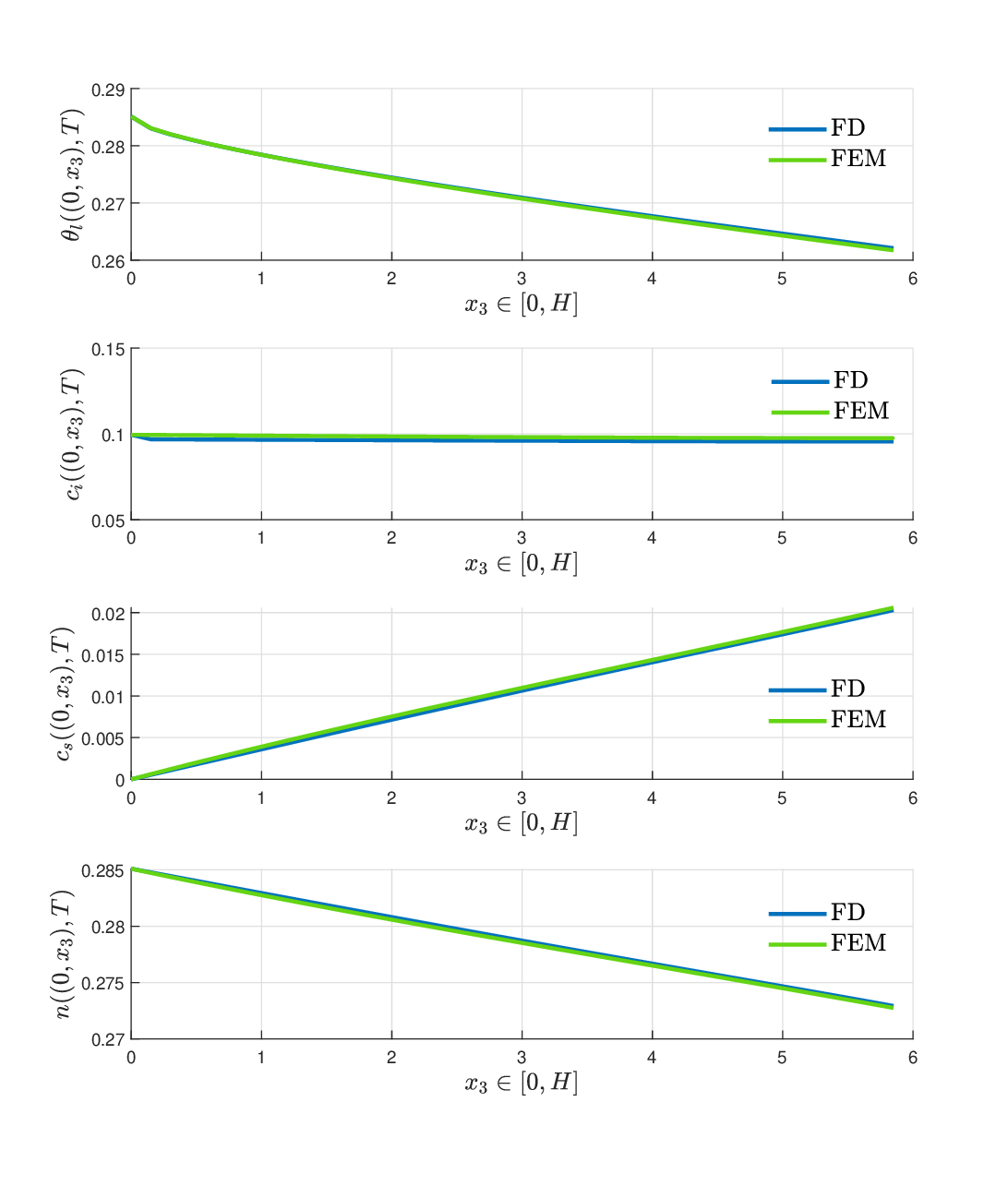}
\vspace{-0.45in}
\caption{Comparison between 2D finite elements (FEM) and 1D finite difference (FD): the approximation of the quantities $\theta_l,\:c_i,\:c_s$ and $n$ is shown, for both the employed techniques, at the final instant of water absorption. FEM results 
have been respectively obtained from the implementation of scheme \eqref{theta k}-\eqref{ci k} and fixing
the vertical variable $x_3$ 
in the range $[0\,\textrm{cm},5.85\,\textrm{cm}]$ determined by the height of the bar.}
\label{fig:FEM_AND_FD_2D}
\end{figure}

We solve the imbibiton problem in dimension $d=2$ setting the following data: the strip $\Omega$ has a base of length $L=0.15\,\textrm{cm}$ and height equal to $H=5.85\,\textrm{cm}$. Absorption of salts has been observed up to time $T\approx 1.6 \cdot 10^6\, \textrm{s}$ (almost $4444\, \textrm{hours}$). The time interval $[0,T]$ has been discretized with step $\Delta t=3.2\,\textrm{s}$, for a number of $5\cdot 10^6$ temporal samplings. 
The plots displayed in Figure \ref{fig:FEM_AND_FD_2D} show the behavior of the quantities $\theta_l,\:c_i,\:c_s$ and $n$ at the final time instant $T$ at the interior points with $x_1=0$ in $\Omega$, so the central vertical axis. Approximations have been obtained by performing the numerical method described in the previous section and discretizing the problem unknowns in space variable with linear finite elements corresponding to a discretization of the boundary of $\Omega$ with segments of length $h=h_{x_3}=0.15 \,\textrm{cm}$ in the vertical direction and $h_{x_1}=0.075\,\textrm{cm}$ in the horizontal one, which induces a conforming triangulation of the interior domain. 

The finite difference scheme \eqref{eq:Num_method_1D} has been performed considering a comparable 1D experiment where the domain $\Omega$ is reduced to the interval $x\in[0\,\textrm{cm},5.85\,\textrm{cm}]$, initial constraints and spatial conditions at the extrema of the interval have been set to replicate the initial and boundary constraints \eqref{initial conditions} and \eqref{Neumann-Robin bc} of the 2D bar. The time step $\Delta t$ has been set equal to the one of the FEM simulation, as well as $\Delta x=h$, 
in order to consider the blue lines in Figure \ref{fig:FEM_AND_FD_2D} as the benchmark results. Let us note that the finite difference discretization is subjected to a CFL condition, leading to a critical choice of the discretization parameters that allows stable 1D results. The plots of the final behavior of the quantities $\theta_l,\:c_i,\:c_s$ and $n$ obtained with 1D finite differences are shown together with the corresponding FEM results in Figure \ref{fig:FEM_AND_FD_2D}.  

\subsection{The 2D bar: experiment for imbibition and drying phases}\label{The 2D bar: esperiment 2 for imbibition and drying phase}


As a second experiment, we solve the imbibition problem by reducing the absorption phase to
$T=4032000\,\textrm{s}$ (almost $1124$ hours) on the same geometry defined for the experiment in subsection \ref{The 2D bar: esperiment 1 for imbibition phase} and for the same initial and boundary conditions. The results for the drying phase have instead been observed up to $\widetilde{T}=2236800\,\textrm{s}$ (almost $620$ hours).

\begin{figure}[h!]
\centering 
\includegraphics[scale=0.62]{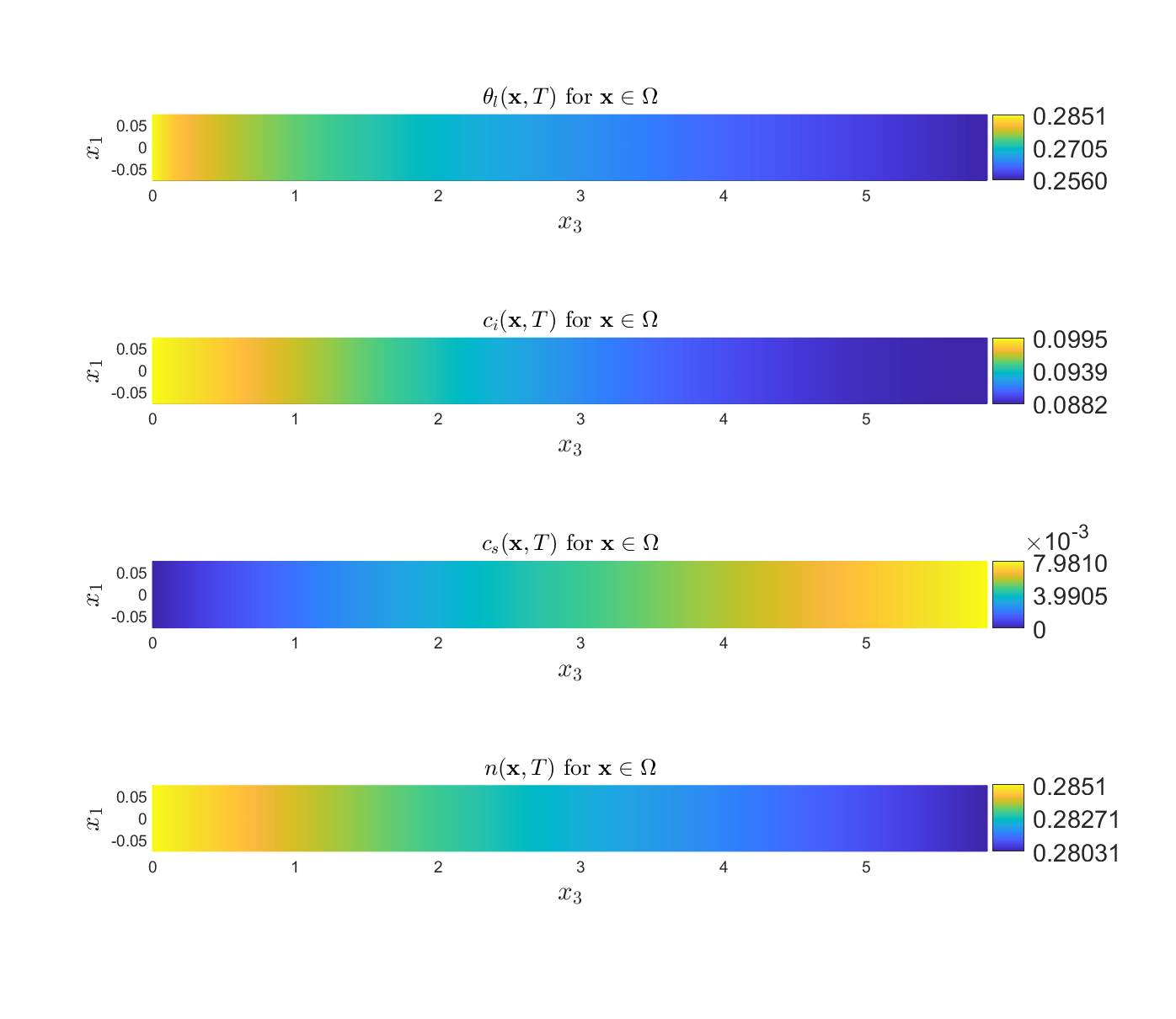}
\vspace{-0.5in}
\caption{\textbf{Water absorption phase.} FEM approximation, at the final instant $T=4032000\,\textrm{s}$, of the fields $\theta_l$, $c_i$, $c_s$ and $n$, all displayed on the whole domain $\Omega=[-0.075,0.075]\times[0,5.85]$.}
\label{fig:imbibition_dt_3_2}
\end{figure}

The discretization parameters are the same for both phases: the strip $\Omega=[-0.15\,\textrm{cm},0.15\,\textrm{cm}]\times[0,5.85\,\textrm{cm}]$ has been discretized at the boundary with segments of length $h=h_{x_3}=0.15 \,\textrm{cm}$ in the vertical direction and $h_{x_1}=0.075\,\textrm{cm}$ (inducing a conforming interior triangulation).
The time intervals of both phases have been discretized with step $\Delta t=3.2\,\textrm{s}$.\\
The colormaps displayed in Figure \ref{fig:imbibition_dt_3_2} show the behavior of the quantities $\theta_l,c_i,c_s$ and $n$ at the final time instant $T$ of the imbibition simulation on the whole rectangular domain $\Omega$. 
It is evident that the considered absorption time (almost $46$ days) is sufficient to completely moisten the bar (values for $\theta_l$ are far from the threshold $\bar{\theta}_l$). The other quantities, related to the crystallization phenomena, are characterized by a slower evolution; nevertheless, if we look at the top of the bar, deposit of salt crystals altered the porosity with a decrease of the $1.6\%$ with respect to the baseline $n_0$. 


\begin{figure}[h!]
\centering
\includegraphics[scale=0.62]{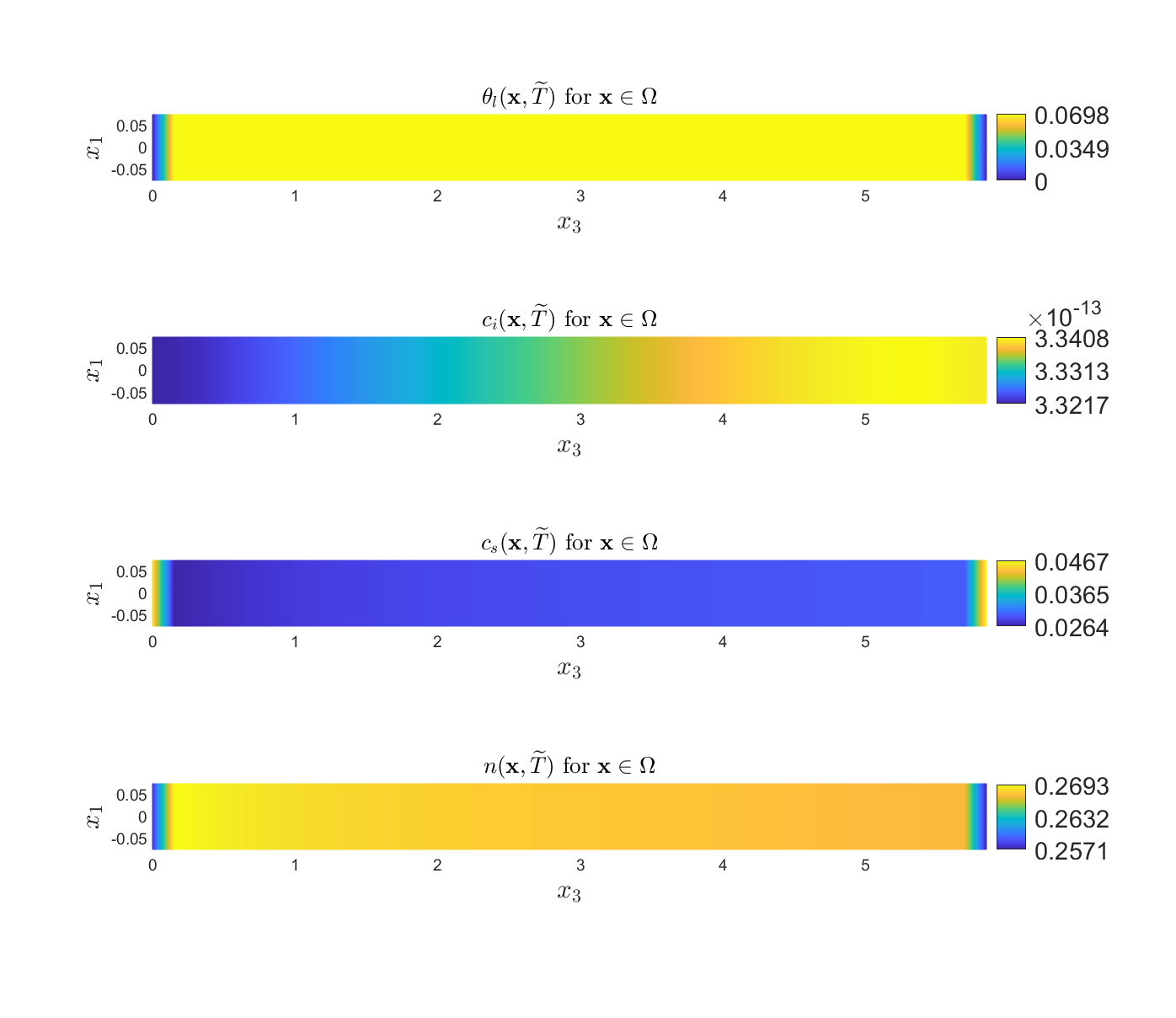}
\vspace{-0.5in}
\caption{\textbf{Drying phase.} FEM approximation, at the final instant $\widetilde{T}=2236800\,\textrm{s}$, of the fields $\theta_l$, $c_i$, $c_s$ and $n$, all displayed on the whole domain $\Omega=[-0.075,0.075]\times[0,5.85]$.}
\label{fig:drying_dt_3_2}
\end{figure}

%
The colormaps shown in Figure \ref{fig:drying_dt_3_2} illustrate the level of the quantities $\theta_l,c_i,c_s$ and $n$ at the final time instant $\widetilde{T}$ of the drying simulation, again on the whole domain $\Omega$. Constraints for $\theta_l$ in \eqref{conditions on theta test 1} in particular induce a reduction of the water volume from the top and the bottom sides of the 2D bar. The free ions concentration $c_i$ at the end of the drying are almost null because of the absence of water supply, while $c_s$ and $n$ are characterized by a maximum magnitude in $\Omega$ of the same order reached at the end of the previous phase of the experiment. As expected, in fact, drying procedure has not significantly changed the quantity of salts already deposited in $\Omega$ and, therefore, the porosity of the material itself.

\subsection{The 3D bar: a realistic experiment of water absorption and drying}

The last experiments consider a realistic geometry in dimension $d=3$. We considered both the absorption phenomenon, observed up to time $T=96000\, \textrm{s}$, and the drying phase, up to time $\widetilde{T}=32000\, \textrm{s}$, for a prism $\Omega$ with base edges of length $L=0.3 \,\textrm{cm}$ and height $H=0.75 \,\textrm{cm}$. The chosen temporal step is $\Delta t=3.2\,\textrm{s}$ for both the absorption and the drying interval, leading to a total of $3\cdot 10^4$ and  $10^4$ instants to study the two phases, respectively. The discretization approach introduced in Section \eqref{sec: Time discretization and iterative weak problem} has been performed approximating the problem unknowns in spatial variable in a finite elements perspective induced by a uniform discretization of the boundary edges with segments of length $h=1.875\cdot 10^{-2} \,\textrm{cm}$. Related uniform triangulation of the external surface and consequent interior discretization with tetrahedral patches lead to an amount of $47396$ piecewise linear basis functions for the spatial approximation of the unknowns $\theta_l$, $c_i$, $c_s$ and $n$.
%

\begin{figure}[h!]
\centering
\captionsetup[subfloat]{labelformat=empty}
\subfloat[$t=320\:s$]{\includegraphics[scale=0.1]{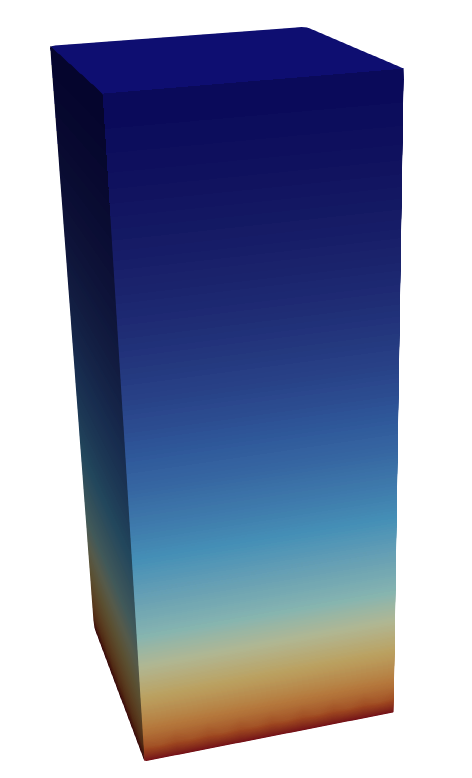}}
\:
\subfloat[$t=640\:s$]{\includegraphics[scale=0.1]{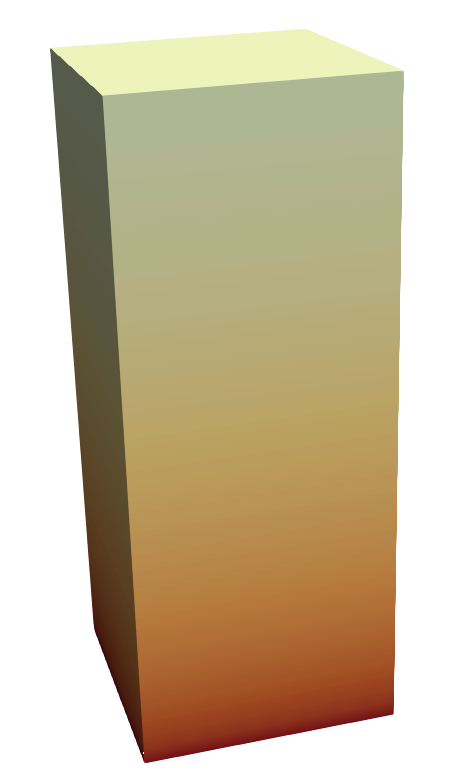}}
\:
\subfloat[$t=960\:s$]{\includegraphics[scale=0.1]{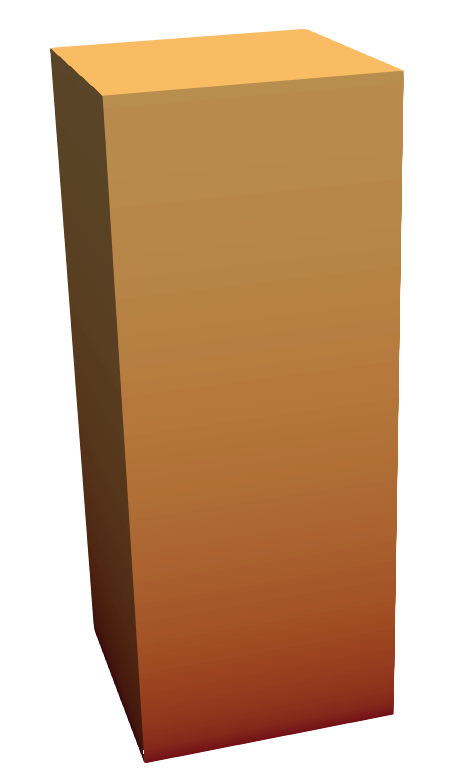}}
\:
\subfloat[$t=1280\:s$]{\includegraphics[scale=0.1]{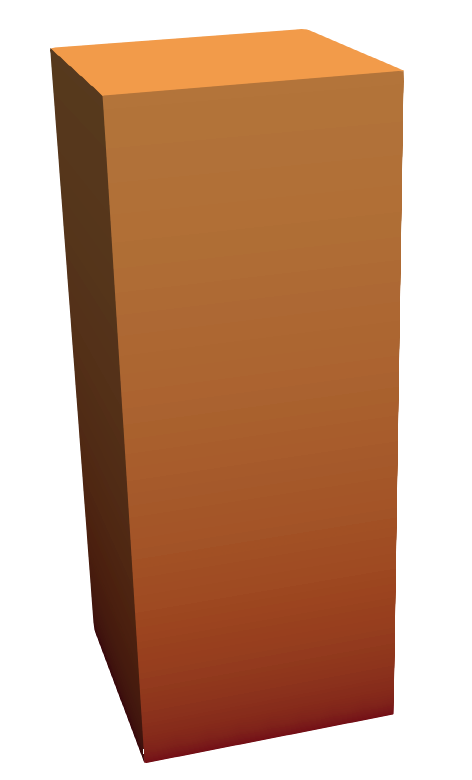}}
\:
\subfloat[$t=1600\:s$]{\includegraphics[scale=0.1]{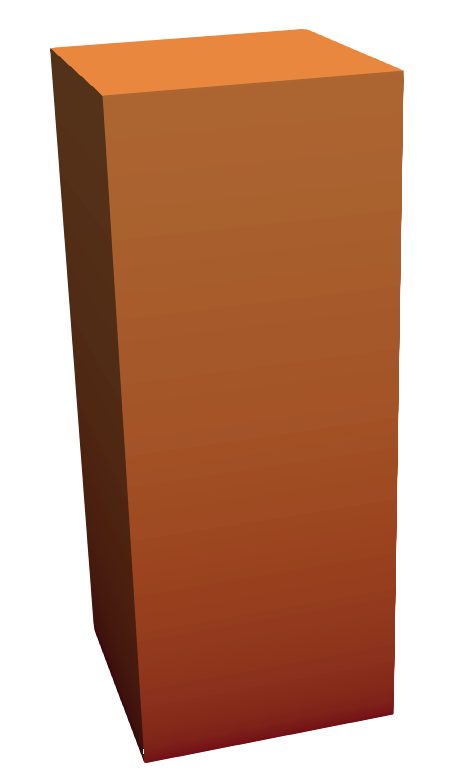}}
\:
\subfloat[$t=1920\:s$]{\includegraphics[scale=0.1]{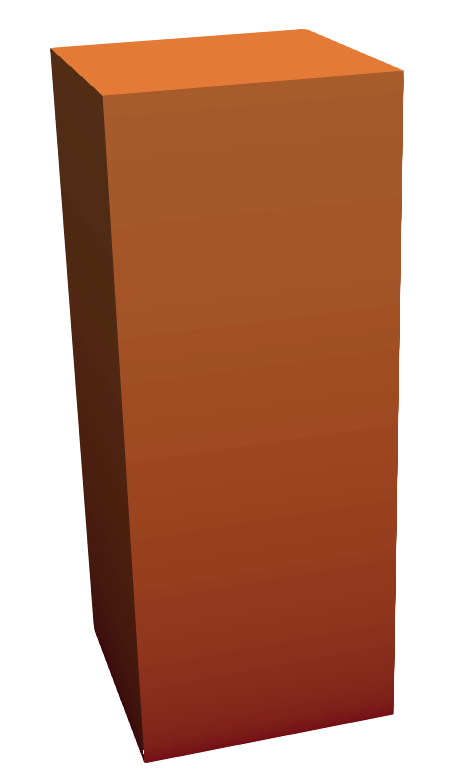}}
\:
\subfloat[$t=96000\:s$]{\includegraphics[scale=0.1]{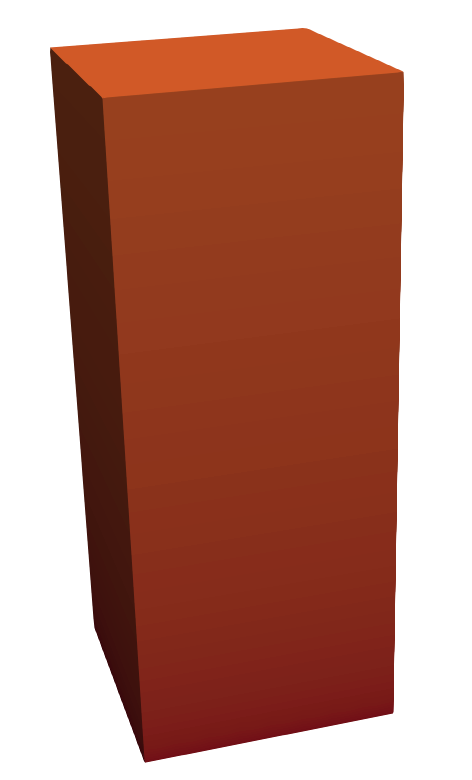}}
\:
\subfloat[]{\includegraphics[scale=0.3]{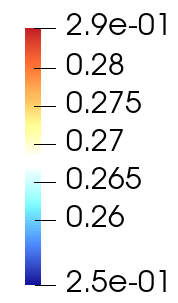}}
\caption{\textbf{Absorption phase}: Approximation by spatial finite elements of water density $\theta_l$ computed on the 3D domain $\Omega$ at different time instants.}
\label{fig:3D experiment theta}
\end{figure}

\begin{figure}[h!]
\centering
\captionsetup[subfloat]{labelformat=empty}
\subfloat[$t=320\:s$]{\includegraphics[scale=0.1]{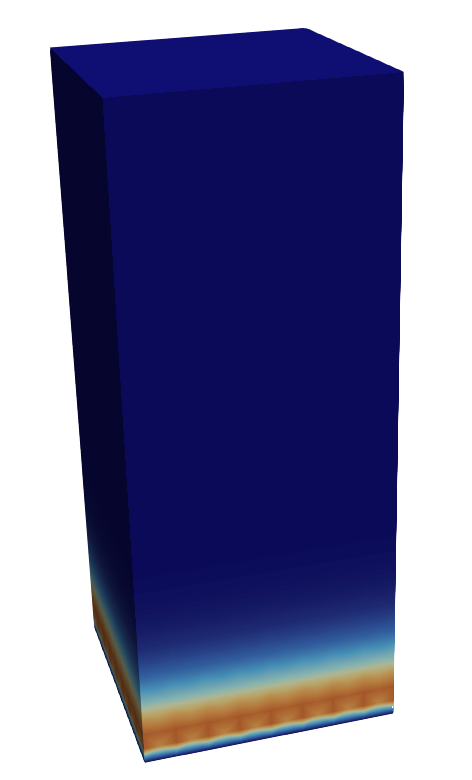}}
\:
\subfloat[$t=640\:s$]{\includegraphics[scale=0.1]{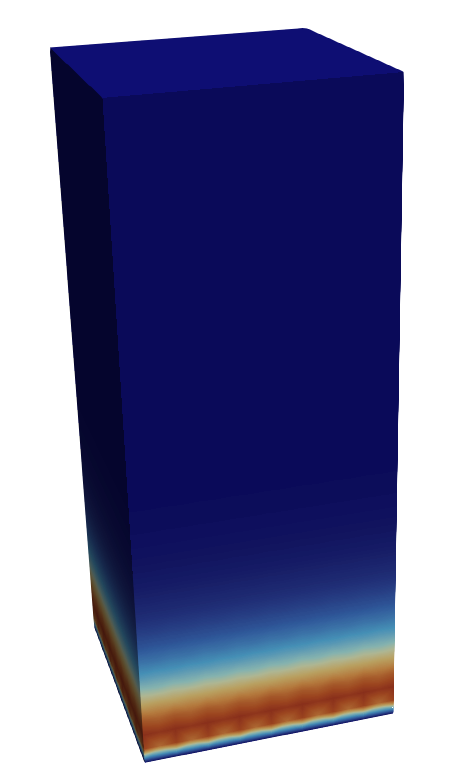}}
\:
\subfloat[$t=960\:s$]{\includegraphics[scale=0.1]{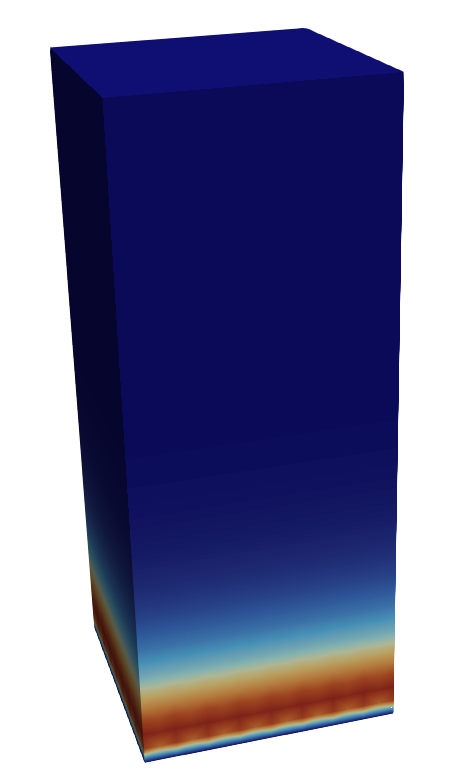}}
\:
\subfloat[$t=1280\:s$]{\includegraphics[scale=0.1]{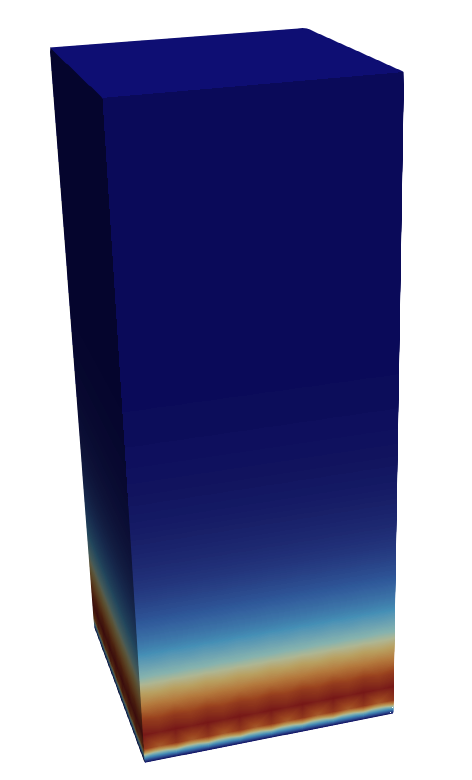}}
\:
\subfloat[$t=1600\:s$]{\includegraphics[scale=0.1]{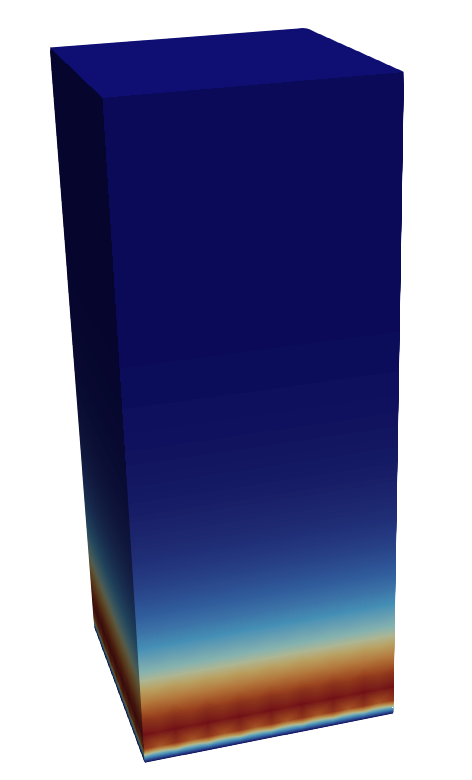}}
\:
\subfloat[$t=1920\:s$]{\includegraphics[scale=0.1]{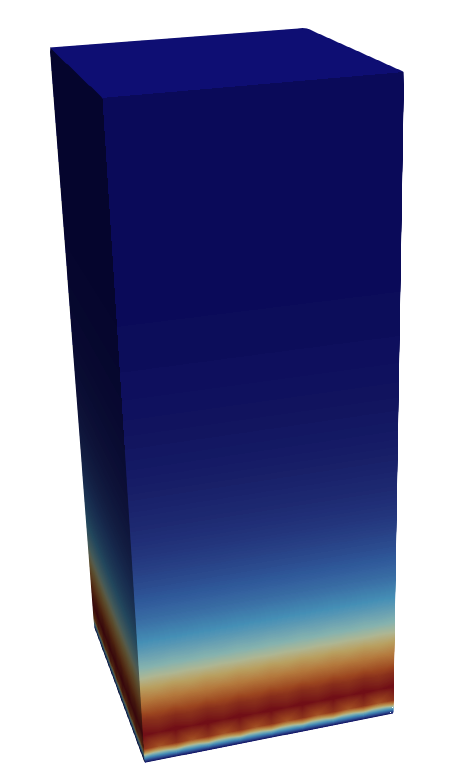}}
\:
\subfloat[$t=96000\:s$]{\includegraphics[scale=0.1]{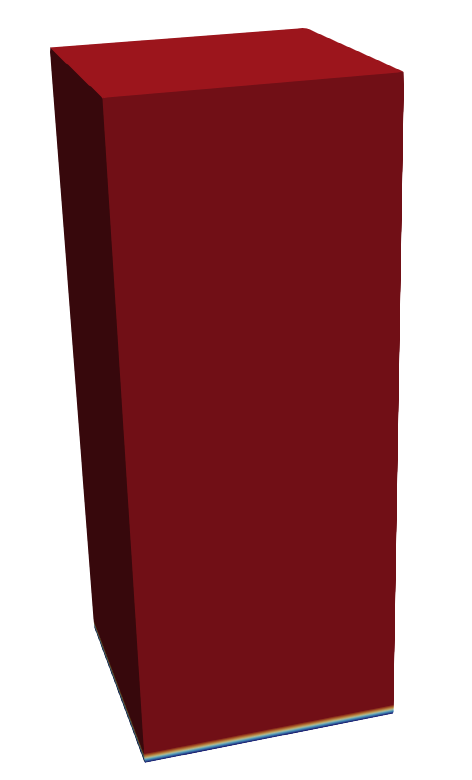}}
\:
\subfloat[]{\includegraphics[scale=0.3]{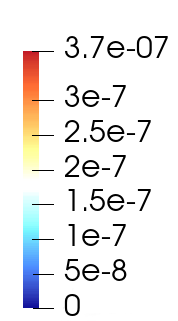}}
\caption{\textbf{Absorption phase}: approximation by spatial finite elements of crystals concentration $c_s$ computed on the 3D domain $\Omega$ at different time instants.} 
\label{fig:3D experiment cs}
\end{figure}

Figures \eqref{fig:3D experiment theta} and \eqref{fig:3D experiment cs} in particular collects approximation results, for the absorption phase, related to $\theta_l$ and $c_s$. Snapshots of their temporal evolution on the surface of the prismatic domain $\Omega$ are displayed for different time instants of simulation. It is observable by the 3D plots that growth of $\theta_l$ in time, namely the fraction of the absorbed water volume, is faster than the deposit of crystal salts through the porous matrix, as visible by the fact that the color of the top extreme for the colorbar in figure \eqref{fig:3D experiment theta} is reached almost starting from $1280\:s$, while, at the same instant, the diffusion front of $c_s$ hasn't reached half of the 3D bar yet.\\


\begin{figure}[h!]
\centering
\captionsetup[subfloat]{labelformat=empty}
\subfloat[$t=320\:s$]{\includegraphics[scale=0.1]{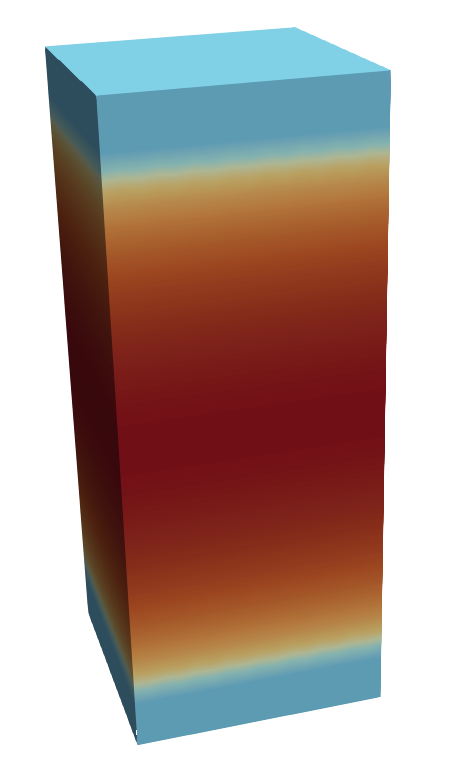}}
\:
\subfloat[$t=640\:s$]{\includegraphics[scale=0.1]{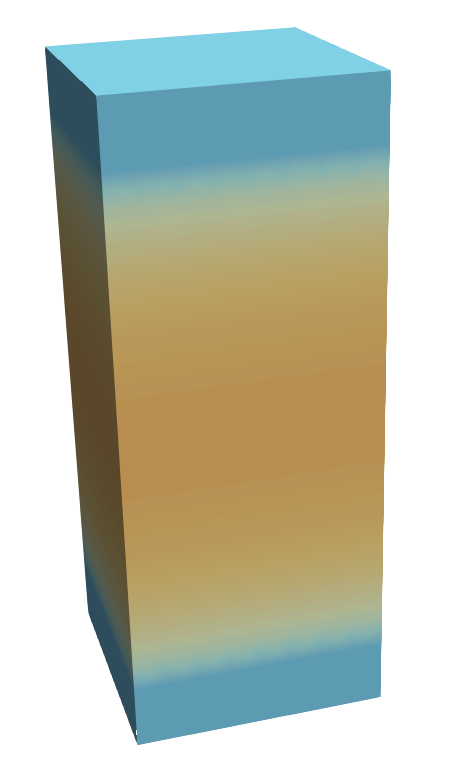}}
\:
\subfloat[$t=960\:s$]{\includegraphics[scale=0.1]{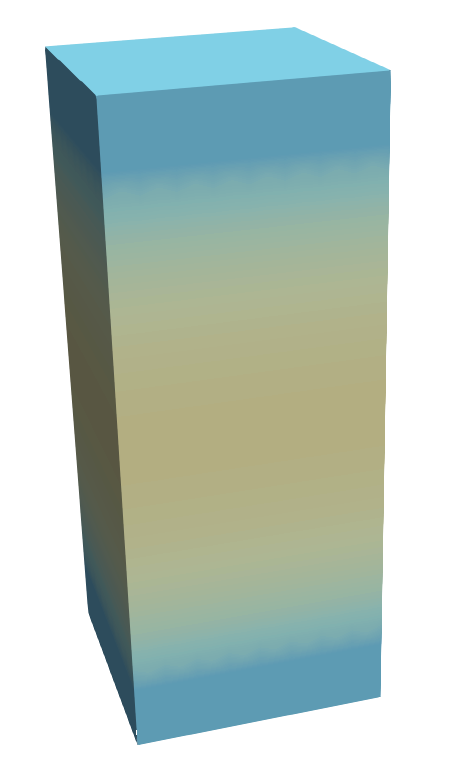}}
\:
\subfloat[$t=1280\:s$]{\includegraphics[scale=0.1]{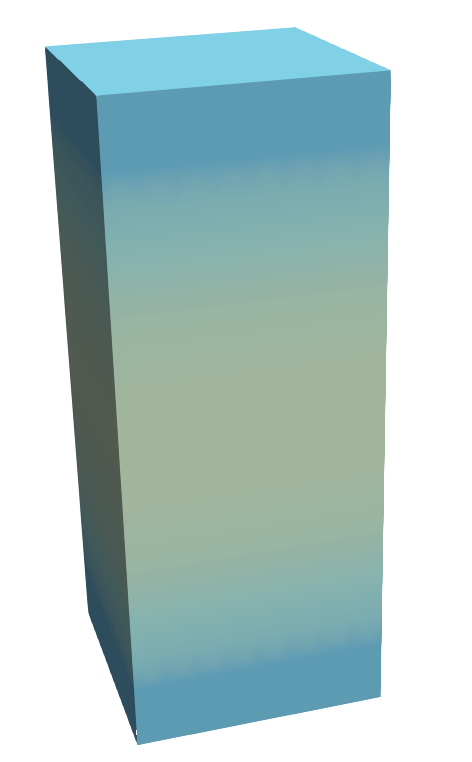}}
\:
\subfloat[$t=1600\:s$]{\includegraphics[scale=0.1]{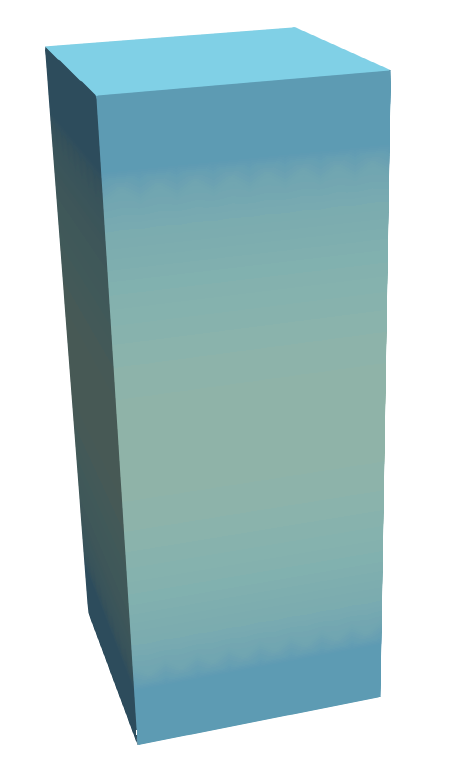}}
\:
\subfloat[$t=1920\:s$]{\includegraphics[scale=0.1]{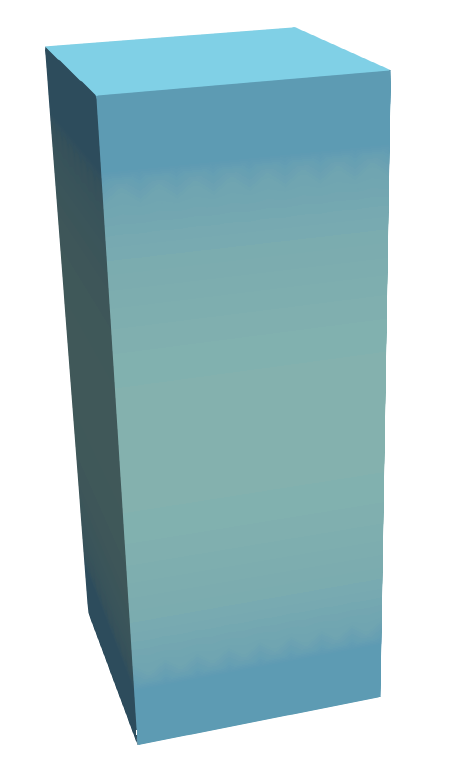}}
\:
\subfloat[$t=32000\:s$]{\includegraphics[scale=0.1]{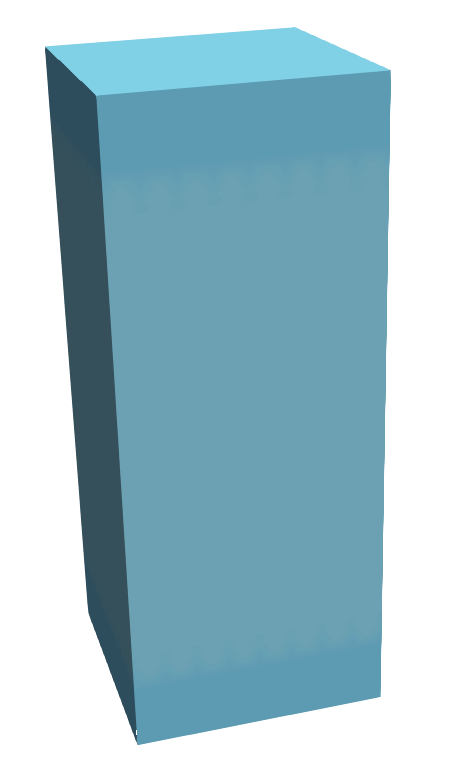}}
\:
\subfloat[]{\includegraphics[scale=0.3]{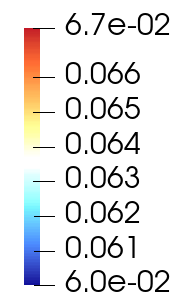}}
\caption{\textbf{Drying phase}: approximation by spatial finite elements of water density $\theta_l$ computed on the 3D domain $\Omega$ at different time instants.}
\label{fig:3D experiment theta dry}
\end{figure}

\begin{figure}[h!]
\centering
\captionsetup[subfloat]{labelformat=empty}
\subfloat[$t=320\:s$]{\includegraphics[scale=0.1]{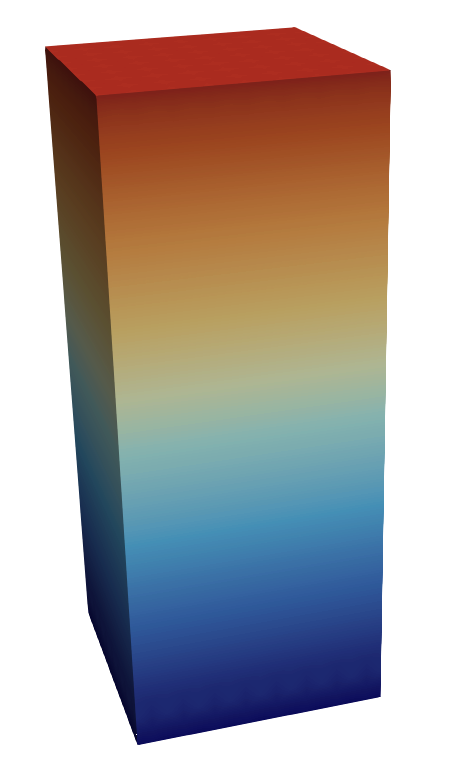}}
\:
\subfloat[$t=640\:s$]{\includegraphics[scale=0.1]{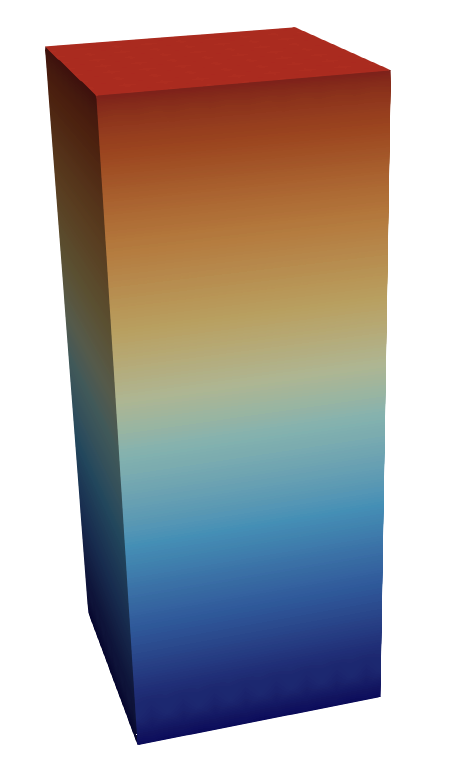}}
\:
\subfloat[$t=960\:s$]{\includegraphics[scale=0.1]{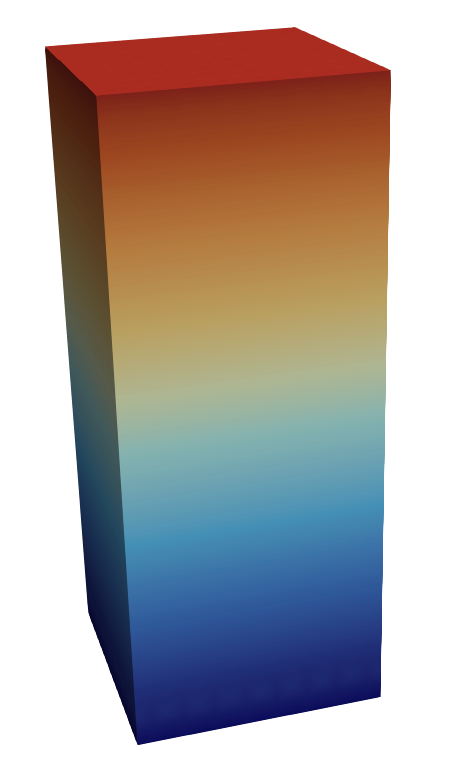}}
\:
\subfloat[$t=1280\:s$]{\includegraphics[scale=0.1]{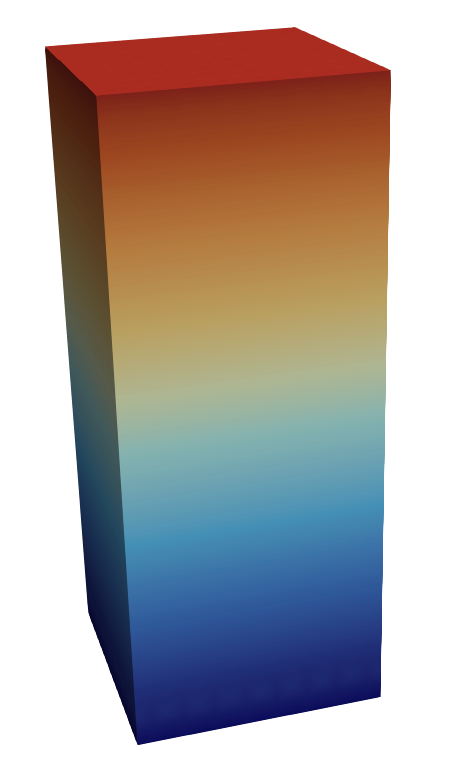}}
\:
\subfloat[$t=1600\:s$]{\includegraphics[scale=0.1]{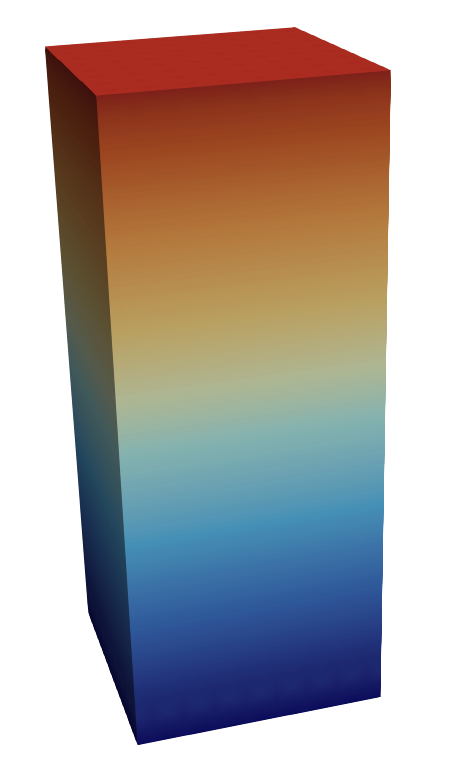}}
\:
\subfloat[$t=1920\:s$]{\includegraphics[scale=0.1]{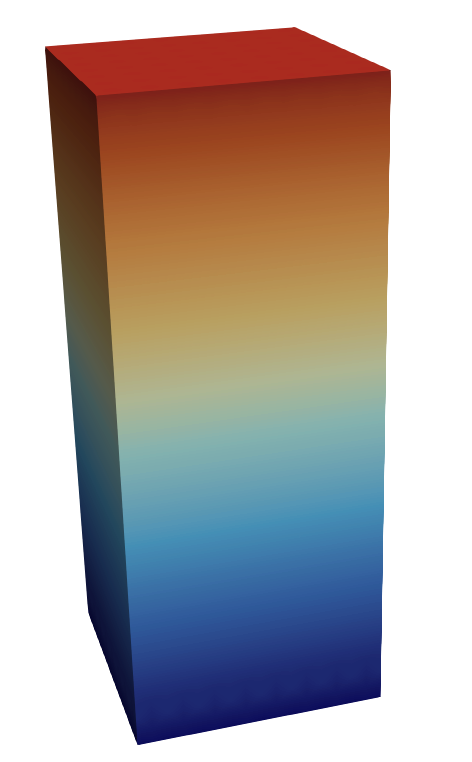}}
\:
\subfloat[$t=32000\:s$]{\includegraphics[scale=0.1]{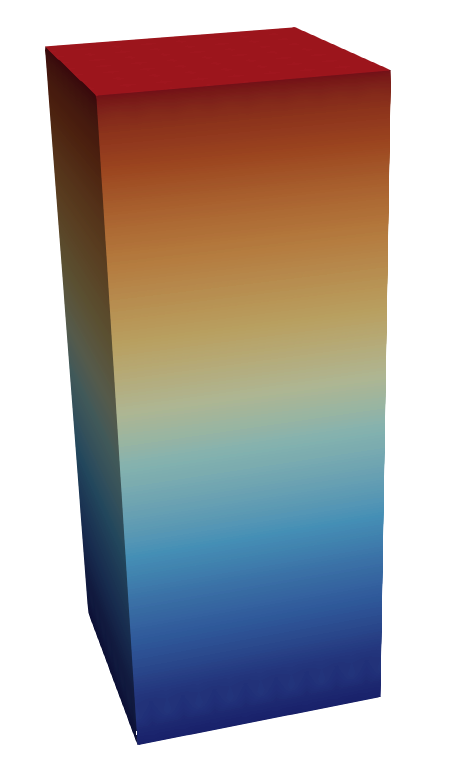}}
\:
\subfloat[]{\includegraphics[scale=0.3]{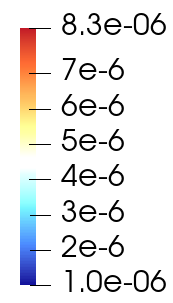}}
\caption{\textbf{Drying phase}: approximation by spatial finite elements of crystals concentration $c_s$ computed on the 3D domain $\Omega$ at different time instants.} 
\label{fig:3D experiment cs dry}
\end{figure}

Snapshots related to the variation during the drying phase, again for $\theta_l$ and $c_s$, are collected in Figures \eqref{fig:3D experiment theta dry} and \eqref{fig:3D experiment cs dry}. Boundary conditions \eqref{conditions on theta test 1} force a decrease of the water volume in $\Omega$ from the bottom and the top basis, allowing in time a rapid drying of the entire 3D prism, as displayed in Figure \ref{fig:3D experiment theta dry}. Similarly to what was observed for the absorption dynamics, the variation of $\theta_l$ is more significant than the other variables. In fact, the frames in Figure \eqref{fig:3D experiment cs dry} show, as expected, that the drying procedure does not alter significantly in time the deposit of crystals obtained at the end of the absorption phase.\\
To quantify the global alteration of porosity at the end of the two phases we chose to employ the metrics for $n$ already introduced in \eqref{eq:Sensitivity_Metrics}, in order to compute the percentage variation of $\mathcal{N}(\gamma,K_s,K_w)$ at the end of the drying phase with respect to its value at beginning of the simulation. In this case, the definition in \eqref{eq:Sensitivity_Metrics} is trivially extended for the 3D geometry considering a volume integral. We observe that $\mathcal{N}=n_0$ at $t=0$ because of the initial conditions for $n$ in \eqref{initial conditions}. In the end, the variation with respect to $n_0$ is small, approximately about $0.00001\%$. This is  justified by the short absorption time, which almost corresponds to one day of water supply. We certainly expect to observe a higher alteration of porosity for a longer experiment.\\

%


\subsection{Experimental convergence analysis}

The final numerical subsection is devoted to the experimental convergence analysis of the numerical method proposed in Section \ref{sec: Time discretization and iterative weak problem}. The convergence orders are expected to be $O(\Delta t)$ in time and, using linear finite elements and $L^2$ norm, $O(h^2)$ in space, this latter for all the unknowns except for the field $\theta_l$, since the related nonlinear PDE can downgrade the quadratic decay of the error (see e.g. \cite{Yang_2015} \textcolor{black}{where, for a nonlinear parabolic PDE, the expected order of convergence in space is proved to be  $O(h^2log^2(2+\frac{1}{h}))\,$}).\\ 
To do that, we solved system \eqref{cont theta}-\eqref{cont ci}, subjected to the absorption conditions \eqref{initial conditions}-\eqref{Neumann-Robin bc}, in the 1D setting. 
The problem has been numerically solved for various levels of space and time discretization, with accuracy tested introducing the following error, in $L^2(\Omega)$ norm at the final instant of simulation,
\begin{equation}\label{error in L2 norm}
    \mathcal{E}(v)=\left(\int_0^H(v(x,T)-v^{\textrm{ref}}(x,T))^2\:dx\right)^{1/2},
\end{equation}
between the unknowns $v=\theta_l,\:c_i,\:c_s$ and $n$ and corresponding reference solutions, in this case numerically computed for a sufficiently small values of the refinement parameters.

\begin{figure}[h!]
\centering
\subfloat[]{\includegraphics[scale=0.5]{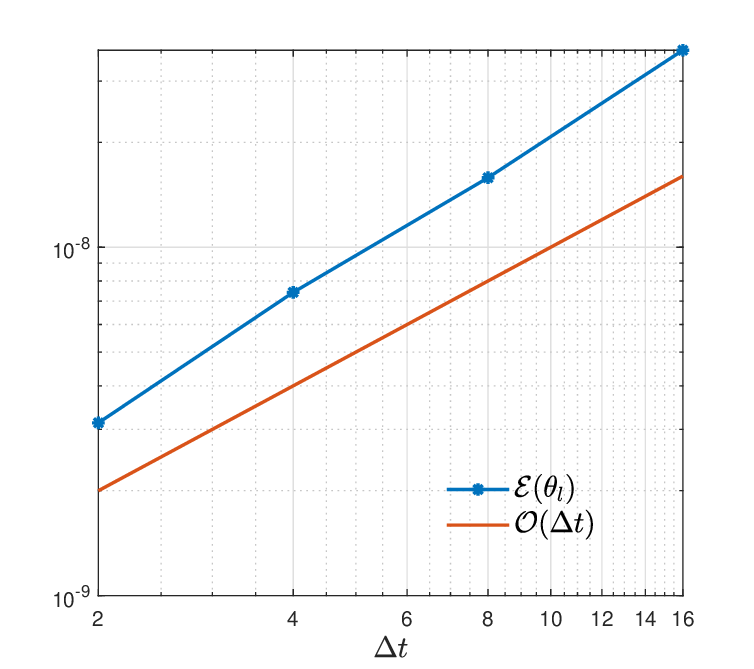}}
\:
\subfloat[]{\includegraphics[scale=0.5]{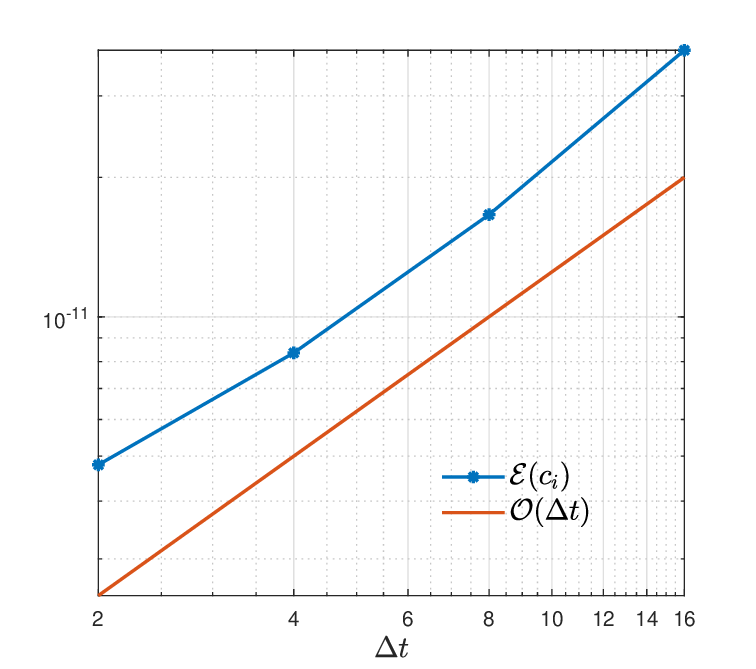}}
\\
\subfloat[]{\includegraphics[scale=0.5]{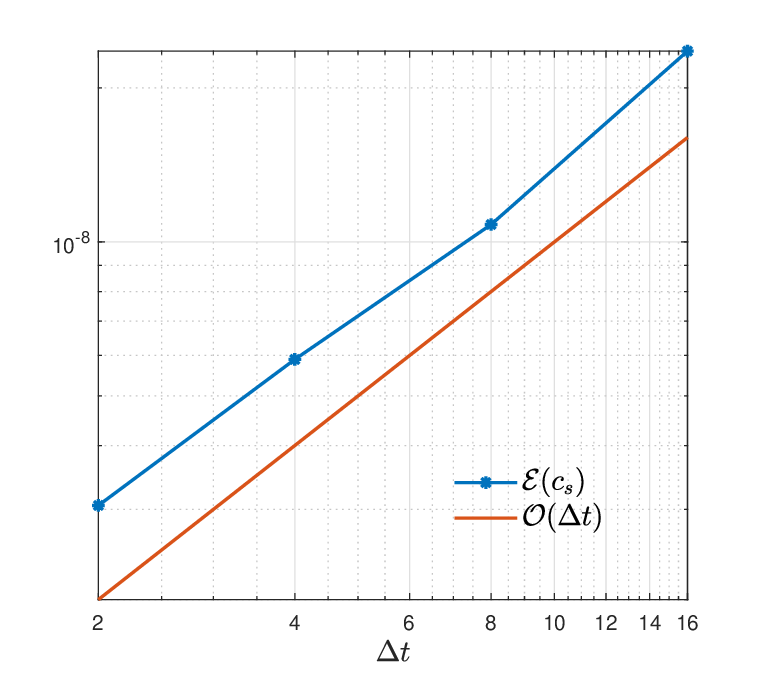}}
\:
\subfloat[]{\includegraphics[scale=0.5]{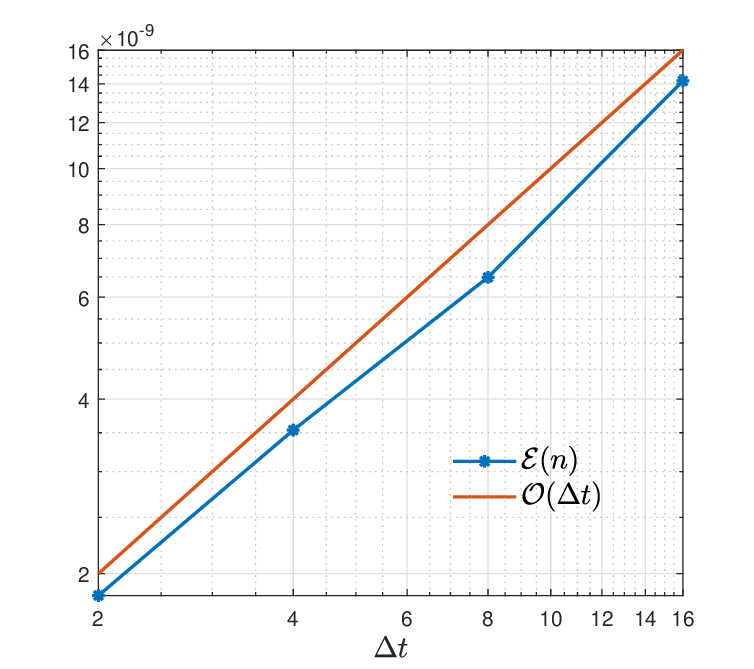}}
\caption{Decay of errors $\mathcal{E}(\theta_l)$ (a), $\mathcal{E}(c_i)$ (b), $\mathcal{E}(c_s)$ (c) and $\mathcal{E}(n)$ (d) considering a fixed space step $h$ and decreasing values of the time step $\Delta t$. Together with the errors, the line $\mathcal{O}(\Delta t)$ is represented.}
\label{fig:convergence error decreasing dt}
\end{figure}

\begin{figure}[h!]
\centering
\subfloat[]{\includegraphics[scale=0.5]
{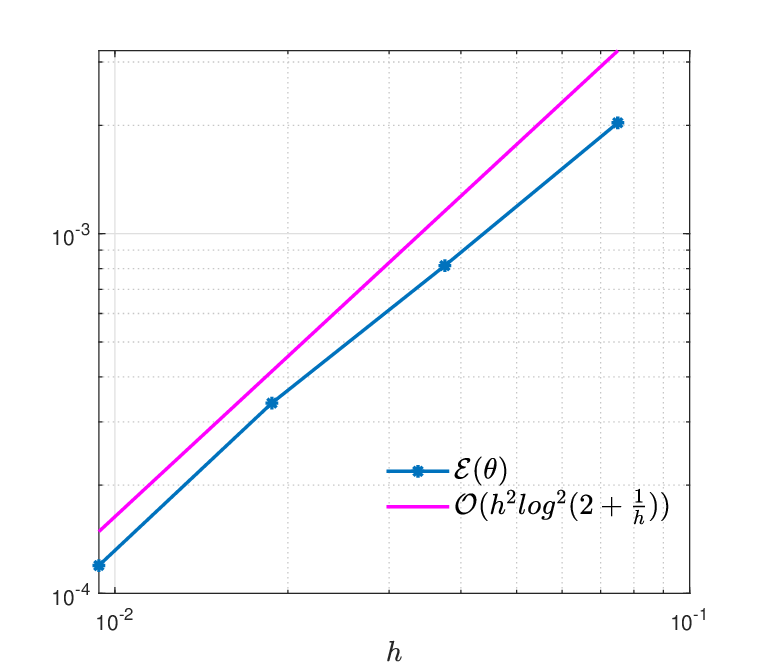}}
\:
\subfloat[]{\includegraphics[scale=0.5]{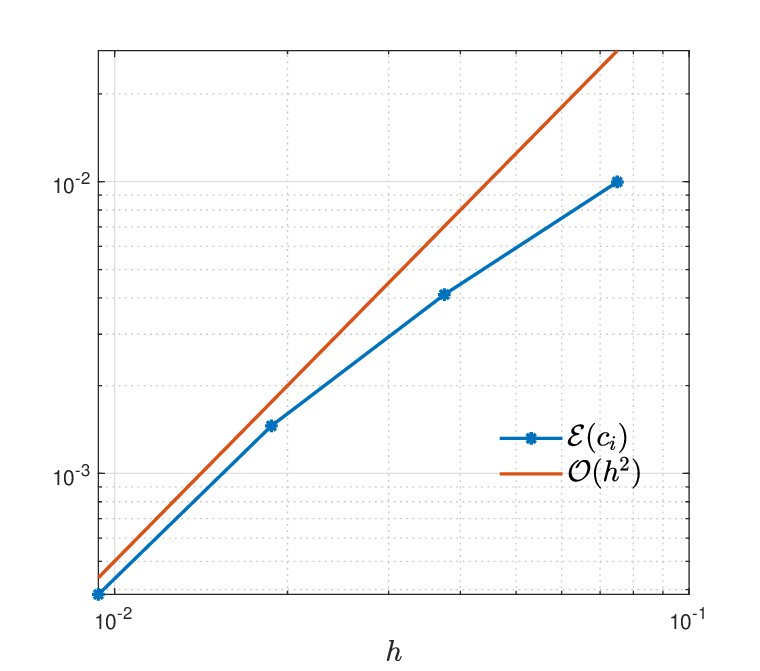}}
\\
\subfloat[]{\includegraphics[scale=0.5]{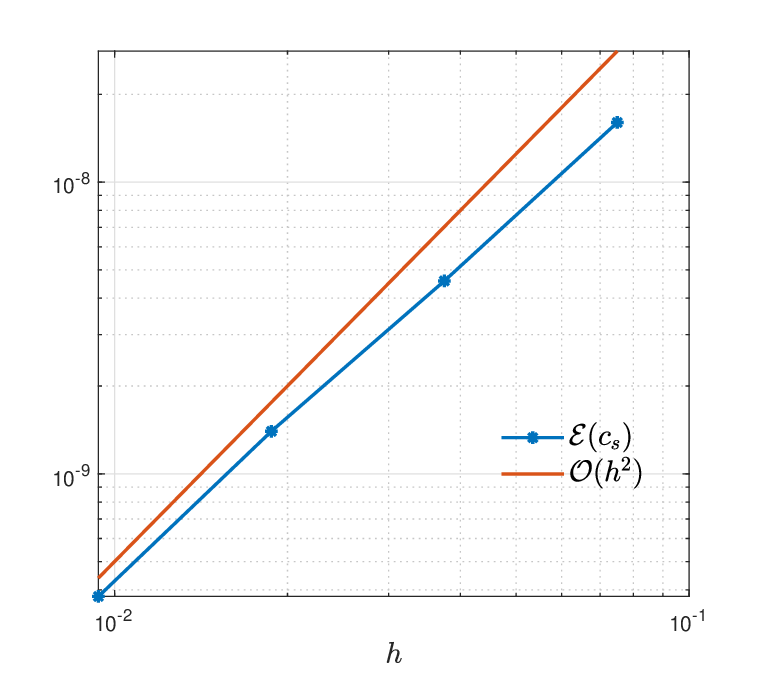}}
\:
\subfloat[]{\includegraphics[scale=0.5]
{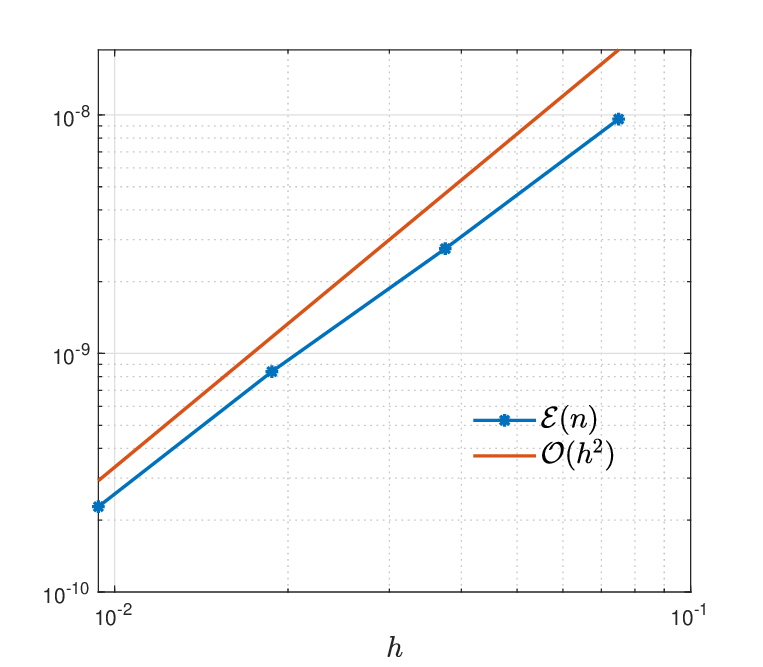}}
\caption{Decay of errors $\mathcal{E}(\theta_l)$ (a), $\mathcal{E}(c_i)$ (b), $\mathcal{E}(c_s)$ (c) and $\mathcal{E}(n)$ (d) considering a fixed time step $\Delta t$ and decreasing values of the space step $h$. \textcolor{black}{Together with the errors, the line $O(h^2)$ is represented in (b), (c), (d), while in (a) a line with slower decay, based on findings in \cite{Yang_2015}, is shown.}}
\label{fig:convergence error decreasing h}
\end{figure}

At first, to recover the convergence order in time, we fix $\Omega$ of length $H=0.6\,\textrm{cm}$ and the final time of analysis as $T=2.56\cdot10^5\,\textrm{s}$, we consider a fixed space refinement of $\Omega$ with $h=3.75\cdot10^{-2}\,\textrm{cm}$ ($16$ space elements) and we implement the implicit time recursion \eqref{theta k}-\eqref{ci k} for $\Delta t=16,8,4,2\,\textrm{s}$. In Figure \ref{fig:convergence error decreasing dt}, the errors $\mathcal{E}(\theta_l)$, $\mathcal{E}(c_i)$ and $\mathcal{E}(c_s)$ are displayed in $log$ scale for each considered time step size, using as reference solutions $\theta_l^{\textrm{ref}}$, $c_i^{\textrm{ref}}$, $c_s^{\textrm{ref}}$ and $n^{\textrm{ref}}$ those obtained using $\Delta t=0.5\,\textrm{s}$. We experimentally deduced a $\mathcal{O}(\Delta t)$ convergence order, as expected by the first order approximation employed in \eqref{discrete time derivatives}.\\
A second test has been run to recover the convergence order in space. Here, we consider $H=0.15\:\textrm{cm}$ and $T=12 \; \textrm{s}$.
The time step has been kept fixed at $\Delta t=$1e-5, a value sufficiently small to not let the analysis be influenced by the time approximation error. The unknowns $\theta_l$, $c_i$, $c_s$ and $n$ have been approximated, in virtue of the implicit time recursion \eqref{theta k}-\eqref{ci k}, for an initial uniform spatial discretization of $\Omega$ with $h=0.075\: \textrm{cm}$ (namely two elements). The procedure has then been repeated halving the space step at each iteration, leading to the plots of $\mathcal{E}(\theta_l)$, $\mathcal{E}(c_i)$, $\mathcal{E}(c_s)$ and $\mathcal{E}(n)$ reported, with respect to the space discretization parameter $h$, in Figure \ref{fig:convergence error decreasing h}. A quadratic decay is experimentally observed for $\mathcal{E}(c_i)$, $\mathcal{E}(c_s)$ and $\mathcal{E}(n)$, justified by the finite element approximation with linear basis functions of the unknowns $c_s$, $n$ and $c_i$ in equations \eqref{cont cs}, \eqref{cont n} and \eqref{cont ci}. For what concerns instead the error $\mathcal{E}(\theta_l)$ in Figure \ref{fig:convergence error decreasing h}(a), in spite of the same approximation of $\theta_l$ field using linear finite elements, its slower convergence is attributable to the non linear behavior of equation \eqref{cont theta}, 
influencing, as commented at the beginning of this section, the  convergence order in space variable for $\theta_l$. 

\begin{remark} From the conducted experimental analysis, we can conjecture that the convergence rate of the FEM approach is superior to the one of the FD scheme recalled in Appendix \ref{app:FD}.
\end{remark}

\section{Conclusions}\label{sec:conclusions}

In this article an extension of the 1D framework in \cite{BRACCIALE201721} is proposed by using a finite-element space discretization with an implicit-explicit in time approximation.  A sensitivity analysis on model parameters calibrated in \cite{BRACCIALE201721} is performed. Moreover, a numerical study for assessing the experimental order of convergence of FEM VS FD method is carried out, see also the Appendix for more details on FD scheme. 
The achievements of the present work can be summarized as follows:
\begin{itemize}
\item a FEM 3D algorithm is developed for the numerical simulation of salt crystallization in porous materials;
    \item  the numerical convergence rate of FEM is conjectured to be higher than FD one for this problem;
    \item the robustness of the model with respect to parameters  exerting the strongest influence on the crystallization dynamics $(\gamma, K_s, K_w)$ is assessed;
    \item the numerical results obtained with FEM algorithm are displayed to show the crystallization process on a brick bar in different experimental settings.
\end{itemize}

 In the future, we will extend the current model to take into account the mechanical damage occurring after a sequence of crystallization cycles, eventually causing important mass decrease of the porous matrix. Moreover, we aim at including the FEM based algorithm on the Stoneverse platform \cite{stoneverse}.

\section*{Acknowledgments}

The authors are grateful to Roberto Natalini for the constructive  discussions on some features of the mathematical model.

This research was supported by PRIN project MATHPROCULT Prot. P20228HZWR.
G.B is responsible of CNR-IAC unit involved in MATHPROCULT project, CUP B53D23015940001.
F.F. is responsible of University of Parma unit involved in MATHPROCULT project, CUP D53D23018820001.
M.P. is supported by the Project PE0000020 CHANGES - CUP B53C22003890006, PNRR Mission 4 Component 2 Investment 1.3, Funded by the European Union - NextGenerationEU under the Italian Ministry of University and Research (MUR). A.A., G. B., G. D. C., C. G. and M. P. are members of the Gruppo Nazionale Calcolo Scientifico - Istituto Nazionale di Alta Matematica (GNCS-INdAM). M.P. is member of the \emph{Research ITalian network on Approximation} (RITA) and of the SIMAI Group \emph{Numerical and Analytical Approximation of Data and Functions with Applications} (ANA$\&$A).

\appendix
\section{Calibration of the top-boundary exchange coefficient} 
\label{app:Kw}

To ensure consistency with the 1D reference model of \cite{BRACCIALE201721}, the exchange coefficient $K_w$ in the Robin boundary condition \eqref{Neumann-Robin bc} is calibrated so that the resulting dynamics of moisture and salt match those obtained with the implicit nonlinear boundary condition (NLBC) used in the reference study.

\begin{figure}[htbp]
  \centering
  \begin{minipage}{0.5\textwidth}
    \centering
    \includegraphics[width=\linewidth]{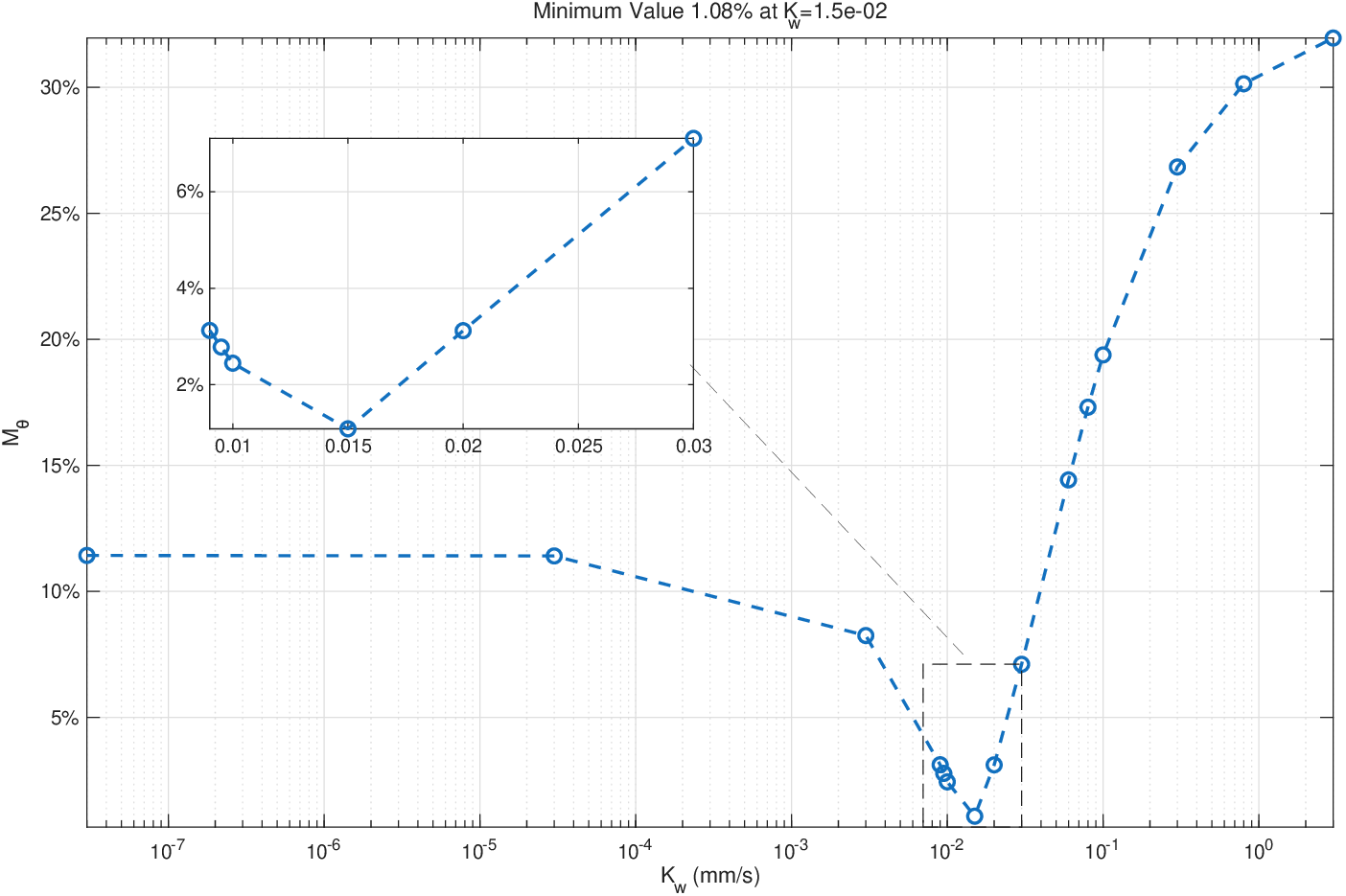}
    \caption{Time-averaged relative error in the total water content for different values of the exchange coefficient $K_w$.}
    \label{fig:Vari_K_min}
  \end{minipage}\hspace{0.05\textwidth}
  \begin{minipage}{0.43\textwidth}
    \centering
    \includegraphics[width=\linewidth]{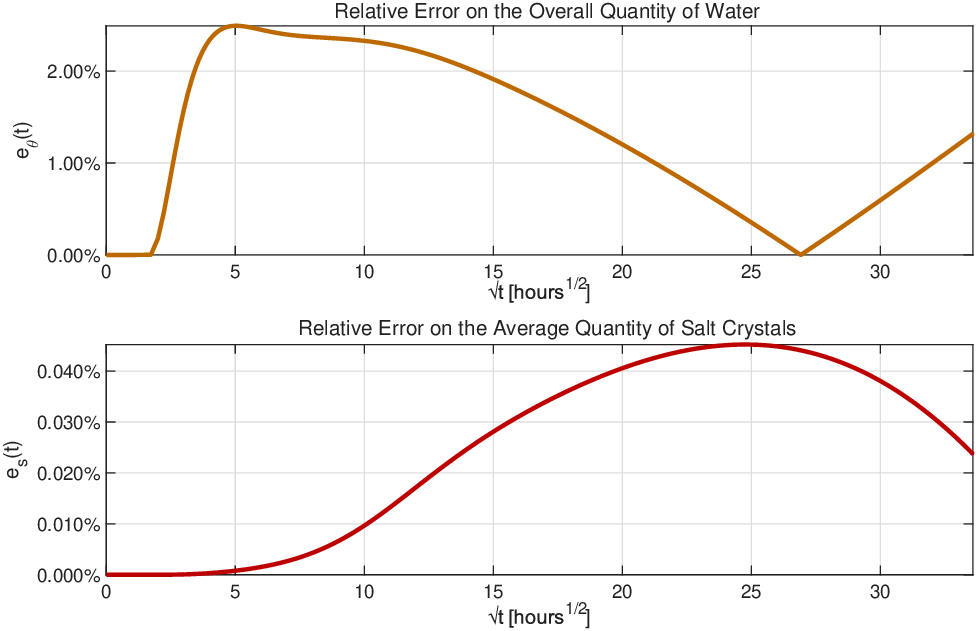}
    \caption{Time evolution of the relative errors for the total water content and salt, comparing the calibrated Robin-based simulation with the reference NLBC solution ($K_w^\ast = 1.5\cdot 10^{-2}$).}
    \label{fig:Intg_diff_Evol}
  \end{minipage}
\end{figure}

Simulations in 1D are performed using the finite-difference scheme described in \eqref{eq:Num_method_1D}. For each choice of $K_w$, the solution is compared with the reference NLBC solution in terms of global quantities, and the relative errors are defined as
\begin{equation}\label{eq:Errori}
    e_\theta(t) =
    \frac{\left|\int_0^H \theta(x,t)\,dx - \int_0^H \theta^{\mathrm{ref}}(x,t)\,dx\right|}{\int_0^H \theta(x,t)\,dx}, 
    \qquad
    e_s(t) =
    \frac{\left|\int_0^H c_s(x,t)\,dx - \int_0^H c_s^{\mathrm{ref}}(x,t)\,dx\right|}{\int_0^H c_s(x,t)\,dx},
\end{equation}
where the superscript $\mathrm{ref}$ denotes the NLBC solution \textcolor{black}{and the subscript for the function $\theta(x,t)$ has been omitted}. Furthermore, the time-averaged discrepancy in total water content is defined as
\begin{equation*}
    M_\theta = \frac{1}{T}\int_0^T e_\theta(t)\,dt.
\end{equation*}
Here, the integral terms are evaluated using the Second Gregory quadrature rule (see, e.g., \cite{Gregory,Gregory2}). 
Figure~\ref{fig:Vari_K_min} shows that $M_\theta$ exhibits a well-defined minimum for $K_w^\ast = 1.5\cdot 10^{-2}$. With this choice, the Robin-based simulation closely reproduces the global dynamics of the reference NLBC throughout the imbibition--drying process, both in terms of moisture content and salt crystallization (see Figure~\ref{fig:Intg_diff_Evol} for the time evolution of the relative errors). This calibrated value is therefore adopted throughout this work.

\section{The finite difference scheme}\label{app:FD}
In this appendix, we provide a numerical method based on FD approximations for the simulation of system \eqref{cont theta}--\eqref{cont ci}. 
The one-dimensional FD framework presented here is adopted to perform sensitivity analysis and investigate boundary conditions in a computationally efficient setting. 
In view of its use within inverse problems and calibration procedures, which require repeated long-time simulations, we resort to an explicit, computationally non-demanding, and easily implementable scheme.
This setting allows us to assess the robustness of the 1D model (see Section \ref{sec:Sensitivity}) and to consistently retain, in the multi-dimensional extension, the parameters calibrated on experimental data in \cite{BRACCIALE201721}. 
Let $\Delta t$ and $\Delta x$ denote positive time and space step sizes, respectively. On the space-time domain $[0,H]\times [0,T]$, consider uniform meshes $t_k=k \Delta t,$ for $k=0,\dots, N_t$ and $x_j=j \Delta x,$ for $j=0,\dots,N_x,$ with $N_t$ and $N_x$ positive integers such that $T=N_t \Delta t$ and $H=N_x \Delta x.$ From now on, for the sake of brevity, we omit the subscript notation for the function $\theta(x,t),$ define the velocity as follows
\begin{equation*}
    V(\theta(x,t),n(x,t))=\left(\dfrac{n(x,t)}{n_0}\right)^2\dfrac{\partial B}{\partial z}\!\left(\dfrac{\theta(x,t)}{n(x,t)}\right),
\end{equation*}
and denote the numerical approximations by 
\begin{equation*}
    \begin{split}
        \theta^{k}_j&\approx \theta(x_j, t_k), \qquad \qquad  c^{k}_{s,j}\approx c_s(x_j, t_k), \qquad \qquad  n^{k}_j\approx n(x_j, t_k), \qquad \qquad  c^{k}_{i,j}\approx c_i(x_j, t_k), \\ 
        V^k_j&\approx V(\theta(x_j, t_k),n(x_j, t_k)), \qquad \qquad   B_j^k\approx B\!\left(\dfrac{\theta(x_j, t_k)}{n(x_j, t_k)}\right) \qquad  \qquad \qquad \ \begin{array}{l}
             j=0,\dots, N_x, \\
             k=0,\dots, N_t.
        \end{array}
    \end{split}
\end{equation*}

We discretize equation \eqref{cont theta} by combining a forward FD discretization of the time derivative with the spatial derivative approximation
\begin{equation}\label{eq:DeltaJ}
    \dfrac{\partial}{\partial z}\left(\! r(x_j,t_k) \, \dfrac{\partial}{\partial z}w(x_j,t_k) \right) \approx \Delta_j(r^k,w^k) := \dfrac{(r_j^k + r^k_{j+1})(w^k_{j+1}-w^k_j) - (r^k_{j-1} + r^k_j)(w^k_j-w^k_{j-1})}{2 \Delta x^2},
\end{equation}
for which we refer to \cite{Driollet_Diele_Natalini, BRACCIALE201721}. Therefore, we set, for each $j=1,\ldots,N_x-1$ and $k=0,\dots,N_t-1,$
\begin{equation*}
     \theta^{k+1}_j= \theta^{k}_j +  \Delta t \ \Delta_j ((n^k/n_0)^2, B^k). 
\end{equation*}
As for the growth of crystals and the porosity variation, a straightforward FD discretization of \eqref{cont cs} and \eqref{cont n} yields
\begin{equation*}
     \dfrac{c_{s,j}^{k+1}-c_{s,j}^k}{\Delta t}=K_s \ c_{i,j}^k (n^k_j-\theta^k_j)^2 + \xbar K (c_{i,j}^k - \bar c)_+ \ \theta^k_j \quad \text{and} \quad n^{k+1}_j = n_0 - \gamma c_{s,j}^{k+1}, \quad \text{for} \ \begin{array}{l}
             j=0,\dots, N_x, \\
             k=0,\dots, N_t-1.
        \end{array}
\end{equation*}
Finally, recalling \eqref{eq:DeltaJ}, equation \eqref{cont ci} is discretized, leading to the following numerical scheme
\begin{equation}\label{eq:Num_method_1D}
    \begin{cases}
        \theta^{k+1}_j= \theta^{k}_j +  \Delta t \ \Delta_j ((n^k/n_0)^2, B^k), \ \ \qquad \qquad \qquad \qquad \qquad \qquad\qquad \qquad\quad \, j=1,\ldots,N_x-1, \\[7pt]
        c_{s,j}^{k+1}= c_{s,j}^k + \Delta t \left[K_s \ c_{i,j}^k (n^k_j-\theta^k_j)^2 + \xbar K (c_{i,j}^k - \bar c)_+ \ \theta^k_j\right], \qquad \qquad \qquad \qquad \quad \, j=0,\ldots,N_x, \\[7pt]
        n^{k+1}_j = n_0 - \gamma c_{s,j}^{k+1}, \qquad \qquad \qquad \qquad \qquad \qquad \qquad  \qquad  \qquad \qquad \qquad \qquad \ \    j=0,\ldots,N_x, \\[7pt]
        c_{i,j}^{k+1}= \dfrac{1}{\theta^{k+1}_j}\left\{\theta^k_j c_{i,j}^{k} + \dfrac{\Delta t}{2 \Delta x} \left(|V^k_{j+1}| c_{i,j+1}^k- 2|V^k_{j}| c_{i,j}^k +|V^k_{j-1}| c_{i,j-1}^k\right) + \Delta t \Delta_j(D\theta^k,c_i^k)\right. \\[7pt]
        \phantom{c_{i,j}^{k+1}=\dfrac{1}{\theta^{k+1}_j}\{ } \left. + \dfrac{\Delta t}{2 \Delta x} \left(V^k_{j+1} c_{i,j+1}^k-V^k_{j-1} c_{i,j-1}^k\right) -  \Delta t \left[K_s \ c_{i,j}^k (n^k_j-\theta^k_j)^2 + \xbar K (c_{i,j}^k - \bar c)_+ \ \theta^k_j\right]\right\},\\[7pt]
        \qquad \qquad \qquad \qquad \qquad  \qquad \qquad \qquad \qquad \qquad  \qquad \qquad \qquad \qquad \qquad \qquad \quad j=1,\ldots,N_x-1, 
    \end{cases}
\end{equation}
where, for $k=0,\dots, N_t-1,$  
\begin{equation}\label{Vdiscr}
 V^k_{0} =0, \qquad \qquad V^k_{j} = \left(B\left(\dfrac{\theta^k_{j+1}}{n^k_{j+1}}\right)-B\left(\dfrac{\theta^k_{j-1}}{n^k_{j-1}}\right)\right)\dfrac{(n_j^k/n_0)^2}{2 \Delta x}, \qquad \qquad \quad  j=1,\ldots,N_x-1.
\end{equation}
We remark that the additional term
\begin{equation}\label{stabilization term}
    \frac{|V^k_{j+1}| c_{i,j+1}^k - 2|V^k_j| c_{i,j}^k + |V^k_{j-1}| c_{i,j-1}^k}{2 \Delta x},
\end{equation}
appearing in the discrete sodium sulfate concentration equation provides an upwind-type stabilization of the nonlinear convective flux $\textrm{div}(c_i(x,t)V(\theta(x,t),n(x,t)))$. 
This contribution does not represent a finite-difference approximation of a differential operator in the continuous model \eqref{cont ci}. 
Rather, it arises from the stabilization of the convective term and corresponds, at the level of the modified equation, to a numerical viscosity with coefficient proportional to $|V|\,\Delta x$, which vanishes as the spatial mesh is refined. 
The inclusion of this term improves the robustness of the explicit scheme in convection-dominated regimes and prevents non-physical oscillations of the concentration field. 
Similar stabilization mechanisms are standard in the numerical approximation of transport problems (see, e.g., \cite{Artificial,LeVeque_2002}).

The initial conditions \eqref{initial conditions} and the boundary conditions \eqref{dirichlet bc}--\eqref{Neumann-Robin bc} are discretized as follows
\begin{equation}\label{eq:initial_conditions_NUMERICAL}
    \theta_j^0=\begin{cases}
        n_0 \ \text{if} \ j=0, \\
        \bar{\theta}_l \ \text{if} \ j>0,
    \end{cases} \qquad
    c_{s, j}^0=0, \qquad
    c_{i, j}^0=\begin{cases}
        \bar{c}_i \ \text{if} \ j=0, \\
        0 \ \text{if} \ j>0,
    \end{cases} \qquad n_j^0=n_0, \qquad  j=0,\dots,N_x, 
\end{equation}
\begin{equation}\label{eq:boundary_conditions_NUMERICAL}
    \theta_0^k=n_0^k, \qquad
    \theta_N^k=\dfrac{4\theta^k_{N_x-1}-\theta^k_{N_x-2}+2\Delta x K_w \bar{\theta}_l}{3+2\Delta x K_w}, \qquad
    c_{i,0}^k=\bar{c}_i, \qquad
    c_{i,N_x}^{k}= \frac{4}{3} c_{i,N_x-1}^{k}- \frac{1}{3}c_{i,N_x-2}^{k},
\end{equation}
for $k=1,\dots,N_t$, obtained from a second order, backward in space finite difference approximation of the derivative operator (we refer to \cite{LeVeque_2007,Quarteroni_2007} for a comprehensive overview on FD schemes). Analogous discretizations are employed for the drying conditions \eqref{initial conditions drying}-\eqref{conditions on theta test 1}.\\
Under standard smoothness assumptions, the finite difference scheme \eqref{eq:Num_method_1D} is expected to be convergent with at least first-order accuracy in both time and space. Numerical investigations, omitted here for the sake of conciseness, support this claim.



\begin{thebibliography}{10}

\bibitem{FEniCS}
M.~Alnæs, J.~Blechta, J.~Hake, A.~Johansson, B.~Kehlet, A.~Logg, C.~Richardson, J.~Ring, M.E. Rognes, and G.N. Wells.
\newblock The {FE}ni{CS} {P}roject {V}ersion 1.5.
\newblock {\em Archive of Numerical Software}, 3(100), 2015.
\newblock \href {https://doi.org/10.11588/ans.2015.100.20553} {\path{doi:10.11588/ans.2015.100.20553}}.

\bibitem{sulfation}
G~G Amoroso and V~Fasina.
\newblock {\em Stone decay and conservation; {A}tmospheric pollution, cleaning; and consolidation}.
\newblock Elsevier Science Pub.,New York, NY, 01 1983.
\newblock URL: \url{https://www.osti.gov/biblio/5767685}.

\bibitem{Driollet_Diele_Natalini}
D.~Aregba-Driollet, F.~Diele, and R.~Natalini.
\newblock {A Mathematical Model for the Sulphur Dioxide Aggression to Calcium Carbonate Stones: Numerical Approximation and Asymptotic Analysis}.
\newblock {\em SIAM Journal on Applied Mathematics}, 64(5):1636--1667, 2004.
\newblock \href {https://doi.org/10.1137/S003613990342829X} {\path{doi:10.1137/S003613990342829X}}.

\bibitem{Freddi_2023}
E.~Bonetti, C.~Cavaterra, F.~Freddi, M.~Grasselli, and R.~Natalini.
\newblock A nonlinear model for marble sulphation including surface rugosity and mechanical damage.
\newblock {\em Nonlinear Analysis: Real World Applications}, 73:103886, 2023.
\newblock \href {https://doi.org/10.1016/j.nonrwa.2023.103886} {\path{doi:10.1016/j.nonrwa.2023.103886}}.

\bibitem{BRACCIALE201721}
M.P. Bracciale, G.~Bretti, A.~Broggi, M.~Ceseri, A.~Marrocchi, R.~Natalini, and C.~Russo.
\newblock Mathematical modelling of experimental data for crystallization inhibitors.
\newblock {\em Applied Mathematical Modelling}, 48:21--38, 2017.
\newblock \href {https://doi.org/10.1016/j.apm.2016.11.026} {\path{doi:10.1016/j.apm.2016.11.026}}.

\bibitem{Data_Informed}
E.~C. Braun, G.~Bretti, M.~Di Fazio, L.~Medeghini, and M.~Pezzella.
\newblock {Data-Informed Mathematical Characterization of Absorption Properties in Artificial and Natural Porous Materials}.
\newblock {\em arXiv}, 2025.
\newblock \href {https://doi.org/10.48550/arXiv.2506.07656} {\path{doi:10.48550/arXiv.2506.07656}}.

\bibitem{Bretti2025}
G.~Bretti and C.~M. Belfiore.
\newblock Mathematical modelling of water absorption properties for historical lime-based mortars.
\newblock {\em GEM - International Journal on Geomathematics}, 16(1):18, Sep 2025.
\newblock \href {https://doi.org/10.1007/s13137-025-00275-2} {\path{doi:10.1007/s13137-025-00275-2}}.

\bibitem{ceseri}
G.~Bretti, M.~Ceseri, R.~Natalini, M.~C. Ciacchella, M.~L. Santarelli, and G.~Tiracorrendo.
\newblock A forecasting model for the porosity variation during the carbonation process.
\newblock {\em Int J Geomath}, 13(13), 2022.
\newblock \href {https://doi.org/10.1007/s13137-022-00204-7} {\path{doi:10.1007/s13137-022-00204-7}}.

\bibitem{goidanich}
G.~Bretti, B.~De Filippo, R.~Natalini, S.~Goidanich, M.~Roveri, and L.~Toniolo.
\newblock Modelling the effects of protective treatments in porous materials.
\newblock In Elena Bonetti, Cecilia Cavaterra, Roberto Natalini, and Margherita Solci, editors, {\em Mathematical Modeling in Cultural Heritage}, pages 73--83. Springer International Publishing, 2021.
\newblock \href {https://doi.org/10.1007/978-3-030-58077-3_5} {\path{doi:10.1007/978-3-030-58077-3_5}}.

\bibitem{volume_mach}
G.~Bretti and R.~Natalini.
\newblock Forecasting damage and consolidation: Mathematical models of reacting flows in porous media.
\newblock In Gabriella Bretti, Cecilia Cavaterra, Margherita Solci, and Michela Spagnuolo, editors, {\em Mathematical Modeling in Cultural Heritage}, pages 187--207. Springer Nature Singapore, 2023.

\bibitem{CdS}
M.~Ceseri, R.~Natalini, and M.~Pezzella.
\newblock {An Integro-Differential Model of Cadmium Yellow Photodegradation}.
\newblock {\em SIAM Journal on Applied Mathematics}, 85(6):2591--2610, 2025.
\newblock \href {https://doi.org/10.1137/24M1709704} {\path{doi:10.1137/24M1709704}}.

\bibitem{charola}
A.E. Charola and E.~Wendler.
\newblock An overview of the water-porous building materials interactions.
\newblock {\em Restoration of Buildings and Monuments}, 21(1):55–65, 2015.
\newblock \href {https://doi.org/10.1515/rbm-2015-0006} {\path{doi:10.1515/rbm-2015-0006}}.

\bibitem{Clarelli_Natalini_Nitsch}
F.~Clarelli, R.~Natalini, C.~Nitsch, and M.~L. Santarelli.
\newblock {\em {A Mathematical Model for Consolidation of Building Stones}}, volume Volume 82 of {\em Series on Advances in Mathematics for Applied Sciences}, pages 232--243.
\newblock WORLD SCIENTIFIC, Sep 2009.
\newblock \href {https://doi.org/10.1142/9789814280303_0021} {\path{doi:10.1142/9789814280303_0021}}.

\bibitem{doehne}
E.~Doehne.
\newblock Salt weathering: a selective review.
\newblock {\em Geological Society, London, Special Publications}, 205(1):51--64, 2002.
\newblock \href {https://doi.org/10.1144/GSL.SP.2002.205.01.05} {\path{doi:10.1144/GSL.SP.2002.205.01.05}}.

\bibitem{grondin}
M.~El-Khoury, F.~Grondin, E.~Roziere, Cortas Rachid, and Chehade~Hage Fadi.
\newblock Chemical degradation vs. creep loading vs. hydration processes in cement-based materials immerged in seawater characterized with a multiscale model.
\newblock {\em Mech Time-Depend Mater}, 29(63), 2025.
\newblock \href {https://doi.org/10.1007/s11043-025-09792-x} {\path{doi:10.1007/s11043-025-09792-x}}.

\bibitem{Freddi_2022}
F.~Freddi and L.~Mingazzi.
\newblock A predictive phase-field approach for cover cracking in corroded concrete elements.
\newblock {\em Theoretical and Applied Fracture Mechanics}, 122:103657, 2022.
\newblock \href {https://doi.org/10.1016/j.tafmec.2022.103657} {\path{doi:10.1016/j.tafmec.2022.103657}}.

\bibitem{guarg}
F.~Guarguaglini.
\newblock A degenerate parabolic model for chemical reactions in porous rocks: existence of solutions.
\newblock {\em MEMOCS}, 13(2):167–200, 2025.
\newblock \href {https://doi.org/10.2140/memocs.2025.13.69} {\path{doi:10.2140/memocs.2025.13.69}}.

\bibitem{guarg2}
F.~Guarguaglini.
\newblock Existence of solutions to a boundary problem for one dimensional degenerate parabolic models for chemical reactions in porous rocks, 2026.

\bibitem{Artificial}
J.L. Guermond and B.~Popov.
\newblock {Invariant Domains and Second-Order Continuous Finite Element Approximation for Scalar Conservation Equations}.
\newblock {\em SIAM Journal on Numerical Analysis}, 55(6):3120--3146, 2017.
\newblock \href {https://doi.org/10.1137/16M1106560} {\path{doi:10.1137/16M1106560}}.

\bibitem{hoke}
G.~D. Hoke and D.~L. Turcotte.
\newblock Weathering and damage.
\newblock {\em Journal of Geophysical Research: Solid Earth}, 107(B10):ECV 1--1--ECV 1--6, 2002.
\newblock \href {https://doi.org/10.1029/2001JB001573} {\path{doi:10.1029/2001JB001573}}.

\bibitem{LeVeque_2002}
Randall~J. LeVeque.
\newblock {\em {Finite Volume Methods for Hyperbolic Problems}}.
\newblock Cambridge University Press, 2002.
\newblock \href {https://doi.org/10.1017/cbo9780511791253} {\path{doi:10.1017/cbo9780511791253}}.

\bibitem{LeVeque_2007}
R.J. LeVeque.
\newblock {\em {Finite Difference Methods for Ordinary and Partial Differential Equations}}.
\newblock Society for Industrial and Applied Mathematics, 2007.
\newblock \href {https://doi.org/10.1137/1.9780898717839} {\path{doi:10.1137/1.9780898717839}}.

\bibitem{Gregory}
E.~Messina, M.~Pezzella, and A.~Vecchio.
\newblock A long-time behavior preserving numerical scheme for age-of-infection epidemic models with heterogeneous mixing.
\newblock {\em Applied Numerical Mathematics}, 200:344--357, 2024.
\newblock New Trends in Approximation Methods and Numerical Analysis (FAATNA20>22).
\newblock \href {https://doi.org/10.1016/j.apnum.2023.04.009} {\path{doi:10.1016/j.apnum.2023.04.009}}.

\bibitem{stoneverse}
E.~Onofri, S.~Bizzarro, S.~Tassa, M.~Czech, and G.~Bretti.
\newblock {Models and Methods in Cultural Heritage. StoneVerse: the open-science platform for cultural heritage}.
\newblock In {\em Digital Heritage 2025}, 2025.

\bibitem{Gregory2}
{Pezzella, M.}
\newblock High order positivity-preserving numerical methods for a non-local photochemical model.
\newblock {\em ESAIM: M2AN}, 59(3):1763--1790, 2025.
\newblock \href {https://doi.org/10.1051/m2an/2025041} {\path{doi:10.1051/m2an/2025041}}.

\bibitem{Quarteroni_2007}
A.~Quarteroni, R.~Sacco, and F.~Saleri.
\newblock {\em {Numerical Mathematics}}.
\newblock Springer New York, 2007.
\newblock \href {https://doi.org/10.1007/b98885} {\path{doi:10.1007/b98885}}.

\bibitem{OAT}
S.~Razavi and H.~V. Gupta.
\newblock {What do we mean by sensitivity analysis? The need for comprehensive characterization of “global” sensitivity in Earth and Environmental systems models}.
\newblock {\em Water Resources Research}, 51(5):3070--3092, 2015.
\newblock \href {https://doi.org/10.1002/2014WR016527} {\path{doi:10.1002/2014WR016527}}.

\bibitem{reale}
R.~Reale, L.~Campanella, M.~P. Sammartino, G.~Visco, G.~Bretti, M.~Ceseri, R.~Natalini, and F.~Notarnicola.
\newblock A mathematical, experimental study on iron rings formation in porous stones.
\newblock {\em Journal of Cultural Heritage}, 38:158--166, 2019.
\newblock \href {https://doi.org/10.1016/j.culher.2019.01.012} {\path{doi:10.1016/j.culher.2019.01.012}}.

\bibitem{siegesmund}
S.~Siegesmund, J.~Menningen, and V.~Shushakova.
\newblock Marble decay: towards a measure of marble degradation based on ultrasonic wave velocities and thermal expansion data.
\newblock {\em Environ Earth Sci}, 80(395), 2021.
\newblock \href {https://doi.org/10.1007/s12665-021-09654-y} {\path{doi:10.1007/s12665-021-09654-y}}.

\bibitem{steiger}
M.~Steiger, K.~Linnow, H.~Juling, G.~G{\"u}lker, A.E. Jarad, S.~Br{\"u}ggerhoff, and D.~Kirchner.
\newblock {Hydration of MgSO4·H2O and Generation of Stress in Porous Materials}.
\newblock {\em Crystal Growth \& Design}, 8(1):336--343, 2008.
\newblock \href {https://doi.org/10.1021/cg060688c} {\path{doi:10.1021/cg060688c}}.

\bibitem{wang}
Q.~Wang, C.~Cheng, E.~Agathokleous, Y.~Liu, X.~Li, and X.~Sheng.
\newblock Enhanced diversity and rock-weathering potential of bacterial communities inhabiting potash trachyte surface beneath mosses and lichens {— A case study in Nanjing, China}.
\newblock {\em Science of The Total Environment}, 785:147357, 2021.
\newblock \href {https://doi.org/10.1016/j.scitotenv.2021.147357} {\path{doi:10.1016/j.scitotenv.2021.147357}}.

\bibitem{warscheid}
Th. Warscheid and J.~Braams.
\newblock Biodeterioration of stone: a review.
\newblock {\em International Biodeterioration \& Biodegradation}, 46(4):343--368, 2000.
\newblock Biodeteriation of Cultural Property 2, Part 2.
\newblock \href {https://doi.org/10.1016/S0964-8305(00)00109-8} {\path{doi:10.1016/S0964-8305(00)00109-8}}.

\bibitem{Yang_2015}
C.~Yang.
\newblock Convergence of a linearized second-order {BDF}-{FEM} for nonlinear parabolic interface problems.
\newblock {\em Computers and Mathematics with Applications}, 70(3):265--281, 2015.
\newblock \href {https://doi.org/10.1016/j.camwa.2015.05.006} {\path{doi:10.1016/j.camwa.2015.05.006}}.

\bibitem{FDM-FEM}
Y.~Zhang, D.~M. Pedroso, and L.~Li.
\newblock {FDM} and {FEM} solutions to linear dynamics of porous media: stabilised, monolithic and fractional schemes.
\newblock {\em International Journal for Numerical Methods in Engineering}, 108(6):614--645, 2016.
\newblock \href {https://doi.org/10.1002/nme.5231} {\path{doi:10.1002/nme.5231}}.

\end{thebibliography}
\end{document}